\documentclass{amsart}
\usepackage{amsmath,amssymb,amscd,amsthm}
\usepackage[matrix,arrow,curve]{xy}

\newcommand{\R}         {\mathbb R}

\newcommand{\bM}        {\mathbf{M}}
\newcommand{\bMM}       {\mathbf{M}/M}
\newcommand{\bT}        {\mathbf{T}}
\newcommand{\bTX}       {\mathbf{T}/X}
\newcommand{\bTZ}       {\mathbf{T}/Z}
\newcommand{\bA}        {\mathbf{A}}
\newcommand{\CC}        {\mathbf{C}}
\newcommand{\DD}        {\mathbf{D}}
\newcommand{\EE}        {\mathbf{E}}
\newcommand{\Sets}      {\mathbf{Sets}}
\newcommand{\Euc}       {\mathbf{Euc}}
\newcommand{\Sch}       {\mathbf{Sch}}

\newcommand{\F}         {\mathcal{F}}

\newcommand{\BB}        {\mathcal B}
\newcommand{\GG}        {\mathcal G}
\newcommand{\UU}        {\mathcal U}
\newcommand{\GGo}       {\GG_{0}}
\newcommand{\GGl}       {\GG_{1}}
\newcommand{\uGG}       {\underline{\GG}}
\newcommand{\uMM}       {\underline{\MM}}
\newcommand{\HHH}       {\mathcal H}
\newcommand{\HHHo}      {\HHH_{0}}
\newcommand{\HHHl}      {\HHH_{1}}
\newcommand{\KK}        {\mathcal K}
\newcommand{\MM}        {\mathcal M}
\newcommand{\TT}        {\mathcal T}

\newcommand{\Geff}      {{\mathcal G}_{\mathrm{eff}}}

\newcommand{\Gtop}      {{\mathcal G}_{\mathrm{top}}}

\newcommand{\Iso}       {\mathrm{Iso}}
\newcommand{\IIso}      {\underline{\mathrm{Iso}}}

\newcommand{\Desc}      {\mathrm{Desc}}

\newcommand{\pr}        {\mathrm{pr}}

\newcommand{\uM}        {\underline{M}}
\newcommand{\uN}        {\underline{N}}
\newcommand{\uf}        {\underline{f}}
\newcommand{\uZ}        {\underline{Z}}
\newcommand{\uC}        {\underline{C}}

\newcommand{\la}        {\leftarrow}
\newcommand{\ra}        {\rightarrow}
\newcommand{\rai}       {\stackrel{\sim}{\rightarrow}}
\newcommand{\arr}[3]    {#1 \stackrel{#2}{\rightarrow} #3}
\newcommand{\lra}       {\longrightarrow}

\newcommand{\intersect} {\cap}
\newcommand{\iso}       {\cong}

\DeclareMathOperator{\Ob}{Ob}
\DeclareMathOperator{\Ar}{Ar}
\DeclareMathOperator{\Hom}{Hom}
\DeclareMathOperator{\HOM}{HOM}

\DeclareMathOperator{\Nat}{Nat}

\DeclareMathOperator{\Op}{Op}
\DeclareMathOperator{\Prin}{Prin}

\DeclareMathOperator{\Maps}{Maps}
\DeclareMathOperator{\Sh}{Sh}
\DeclareMathOperator{\St}{St}
\DeclareMathOperator{\EG}{EG}
\DeclareMathOperator{\PEG}{PEG}

{\swapnumbers
  {\theoremstyle{plain}
  \newtheorem{theorem}{Theorem}
  \newtheorem{Def}[theorem]{Definition}
  \newtheorem{DefProp}[theorem]{Definition-Proposition}
  \newtheorem{proposition}[theorem]{Proposition}
  \newtheorem{cor}[theorem]{Corollary}
  \newtheorem{lemma}[theorem]{Lemma}
  }
  {\theoremstyle{remark}
    \newtheorem{remark}[theorem]{Remark}
    \newtheorem{example}[theorem]{Example}
    \newtheorem*{note}{Note}
  }
}

\begin{document}

\title{Topological and Smooth Stacks}

\author{David S. Metzler}

\begin{abstract}
We review the basic definition of a stack and apply it to the
topological and smooth settings. We then address two subtleties
of the theory: the correct definition of a ``stack over a stack''
and the distinction between small stacks (which are algebraic objects)
and large stacks (which are generalized spaces). 
\end{abstract}

\maketitle

\section{Introduction}
\label{sec:intro}

This paper has two purposes. The first is to present the basic definitions
related to stacks in the topological and smooth categories. 
These notions are standard in the setting of algebraic geometry
but are less well-known in the smooth and topological
settings. Second, we want to clarify some subtleties in the
standard notions. In particular we discuss the notion of a stack
over another stack and the distinction between
large and small stacks.

We will try to be as concrete as is reasonably possible,
concentrating on the topological and smooth situations.
We will make the basic definitions in the general
setting of arbitrary sites, since in setting up the framework,
it is simpler to do it in this degree of generality.
However we will not go at all deeply into the general theory,
so a cursory reading of the general definitions, supplemented
by a thorough understanding of the key examples, should suffice
for the reader most interested in the topological and smooth applications.
Also, when discussing the topological and smooth cases, we will 
do both in parallel where that is possible.

We will make some comments on one of the most interesting special
cases, namely orbifolds, which, properly understood, are
the smooth analogue of Deligne-Mumford stacks.
The treatment here is intended to be complementary to e.g.
\cite{Moerdijk:OrbGrpIntro} and the recent preprint
\cite{TuXuLaurent:TwistedKStacks}. We leave for a later paper
a combination of the two issues presented here, namely
the notion of a small stack over a large stack. 

\vspace{5mm}

The motivating philosophy of stacks is simple. Let us discuss it
in the case of smooth manifolds for concreteness. We start with the
observation that any manifold $M$ is determined up to canonical isomorphism
if we know all smooth maps into $M$ from any other manifold $N$.
In fact, since manifolds are locally Euclidean, we know $M$ if
we know all smooth maps from $\R^{n}$ to $M$ for all $n$. (For example,
the set of smooth maps from $\R^{0} = \ast$ to $M$ is just the set
of points of $M$. The set of smooth 
maps from $R^{n}$ to $M$, for $n = \dim M$,
includes all of the local charts for $M$. This data clearly determines
$M$ up to diffeomorphism.) This is Grothendieck's philosophy of
the ``functor of points:'' he considers a map $N \ra M$ as
a generalized ``point'' of $M$. The \textit{representable functor}
$\uM$ from the category of manifolds to the category of sets
given by $\uM(N) = \{f:N \ra M \}$ and 
$\uM(N' \stackrel{g}{\ra} N) = (f \mapsto f \circ g)$ considered
as a substitute for $M$ itself. One is then led to consider
other functors from (Manifolds) to (Sets) which are not necessarily
representable, as some kind of a generalized space. What we will discover
is that we need to generalize a bit further to get a nice description
of, for example, orbifolds.

To be complete, let us state the fundamental (but trivial) lemma
underlying this idea, the \textit{Yoneda Lemma}.
Let $\CC$ be a category and let $C$ be an object of $\CC$.
We denote by $\uC$ be the contravariant functor represented by $C$, 
given by 
\begin{align*}
  \uC:             \CC       &\ra      \Sets                  \\
                    A        &\mapsto  \{f: A \ra C\}         \\
     (A \stackrel{g}{\ra} B) &\mapsto  (f \mapsto f \circ g).
\end{align*}

\begin{lemma}\label{lem:yoneda}
  Let $\CC$ be a category, let $F: \CC \ra \Sets$ be a contravariant
  functor, and let $C$ be an object of $C$. Then there is a natural
  bijection
\[
   \phi: F(C) \ra \Nat(\uC,F)
\]
  given by $\phi(x)(A)(\alpha) = F(\alpha)(x),$ where
  $x \in F(C)$ and $\alpha: A \ra C$. 
  Here $\Nat(\uC,F)$ is the set of natural transformations.
  The inverse is given by $\phi^{-1}(\psi) = \psi_{C}(1_{C})$.
\end{lemma}

The following corollary is a little more intuitive. It is
just the lemma applied to another representable functor
$F = \uC'$.
\begin{cor}\label{cor:yoneda}
  Let $\CC$ be a category and let $C, C'$ be objects of $\CC$.
  Let $\uC, \uC'$ be the functors represented by $C, C'$ respectively
  Then there is a natural bijection 
  $\phi$ from $\Hom_{\CC}(C,C')$ to the set of natural transformations 
  of functors $\Nat(\uC,\uC')$, given by 
  \[
    \phi(C \stackrel{h}{\ra} C')(A) = (f \mapsto h \circ f), 
      \quad \text{~where~} f:A \ra C. 
  \]
  The inverse of $\phi$ is given by evaluation at the identity $1_{C'}$.
\end{cor}

We note one convention: we will often say $C \in \CC$ to mean $C \in \Ob(\CC)$
where there is likely to be no confusion.

\section{Stacks}
\label{sec:stacks}

We begin with the general notion of a \textit{site}, which
is the correct general setting for sheaves and stacks.
We will first present the definition of a site based on sieves,
as this often provides the most elegant constructions.
The following three definitions, originally due to Grothendieck 
\cite{GT}, \cite{SGA4.1},
are taken, with minor modification, from MacLane and Moerdijk \cite{MacMoer:Sheaves}.

\begin{Def}\label{def:sieve}
  Given a category $\CC$ and an object $C \in \Ob \CC$,
  a \textbf{sieve} (French ``crible'') $S$ on $C$ is a family of arrows
  of $\CC$, all with target $C$, such that 
  \[
    f \in S \implies fg \in S
  \]
  whenever $fg$ is defined. (I.e. $S$ is a right ideal under composition.)
  Given a sieve $S$ on $C$ and an arrow $h: D \ra C$, we define
  the \textbf{pullback} sieve $h^{*}(S)$ by 
  \[
    h^{*}(S) = \{g \: | \: \mathrm{target}(g) = D, \: hg \in S \}. 
  \]
\end{Def}

\begin{Def}\label{def:site}
  A \textbf{site} $(\CC,J)$ is a small category $\CC$
  equipped with a \textbf{Grothendieck topology} $J$,
  that is, a function $J$ which assigns to each object
  $C$ of $\CC$ a collection $J(C)$ of sieves on $C$, called 
  \textbf{covering sieves}, such that 
  \begin{enumerate}
  \item the maximal sieve $t_{C} = \{ f \: | \: \mathrm{target}(f) = C \}$
        is in $J(C)$;
  \item (stability) if $S \in J(C)$, then $h^{*}(S) \in J(D)$ for any
        arrow $h: D \ra C$;
  \item (transitivity) if $S \in J(C)$ and $R$ is any sieve on $C$
        such that $h^{*}(R) \in J(D)$ for all $h:D \ra C$ in $S$, 
        then $R \in J(C)$.
  \end{enumerate}
\end{Def}

It is useful to note two simple consequences of these axioms.
First, there is a somewhat more intuitive transitivity property:
\begin{itemize}
\item [$3'$]($\mathrm{transitivity}'$) If $S \in J(C)$ is a covering sieve 
      and for each $f: D_{f} \ra C$ in $S$ there is a covering sieve 
      $R_{f} \in J(D_{f})$, then the set of all composites $f \circ g$,
      where $f \in S$ and $g \in R_{f}$, is a covering sieve of $C$.
\end{itemize}
Next we have the fact that any two covering sieves have a common refinement,
in fact, their intersection.
\begin{itemize}
\item [4] (refinement) If $R,S \in J(C)$ then $R \intersect S \in J(C)$.
\end{itemize}

It is often more intuitive to work with a \textit{basis}
for a topology (also called a \textit{pretopology}).
\begin{Def}\label{def:basis}
  A \textbf{basis} for a Grothendieck topology on a category
  $\CC$ is a function $K$ which assigns to every object $C$ of $\CC$
  a collection $K(C)$ of families of arrows with target $C$, 
  called \textbf{covering families}, such that
  \begin{enumerate}
\item if $f: C' \ra C$ is an isomorphism, then $\{f\}$ is a covering family;
\item (stability) if $\{f_{i}: C_{i} \ra C \}$ is a covering family, then for
      any arrow $g: D \ra C$, the pullbacks $C_{i} \times D$ exist and 
      the family of pullbacks $\pi_{2}: C_{i} \times D \ra D$ 
      is a covering family (of $D$);
\item (transitivity) if $\{f_{i}: C_{i} \ra C \: |\: i \in I\}$ 
      is a covering family and for each $i \in I$, one has a 
      covering family $\{g_{ij}: D_{ij} \ra C_{i} \: | \: j \in I_{i}\}$,
      then the family of composites 
      $\{f_{i} g_{ij}: D_{ij} \ra C \: | \: i \in I, j \in I_{i}\}$ 
      is a covering family.
\end{enumerate}
\end{Def}

Any basis $K$ generates a topology $J$ by 
\[
  S \in J(C) \Leftrightarrow \exists R \in K(C) \text{~with~} R \subset S.
\]

In other words, the covering sieves on $C$ are those which
refine some covering family $R$. See Example \ref{ex:smalltopsite}
below. Usually we will describe sites in terms of a basis.
We will often abuse notation and refer to a site $(\CC,J)$ simply
as $\CC$.

One simple way in which new sites arise is the \textit{induced site}.
\begin{Def}\label{def:inducedsite}
  Let $(\CC,J)$ be a site and let $u: \bA \ra \CC$ be a functor.
  Assume that $u$ preserves all pullbacks that exist in $\bA$.
  The \textbf{induced topology} $J|_{A}$ on $\bA$ is defined
  in terms of the following basis. A family $f_{i}:A_{i} \ra A$
  is a covering family for the induced topology if and only if 
  the family $u(f_{i}):u(A_{i}) \ra u(A)$ is a covering family
  for $J$.
\end{Def}
\begin{example}\label{ex:inducedsubcategory}
  Let $(\CC,J)$ be a site and let $\bA \subset \CC$ be a full
  subcategory. Assume the inclusion functor preserves all pullbacks
  that exist in $\bA$.
  Then the induced topology on $\bA$ will
  also be called the \textit{restriction} of $J$ to $\bA$
  and will be denoted $J|_{A}$.
\end{example}

We now present key examples of sites.
Let us first get out of the way a set-theoretic issue.
We will want to discuss, for example, ``the category of stacks on 
the category of all topological spaces,'' but strictly
speaking this does not exist, since the category
of topological spaces does not have a set of objects,
but rather a proper class. To avoid this problem we will
consider throughout some fixed category $\bT$ of topological
spaces which has a set of objects, or at least, is equivalent
to such a category. 

For example, the category of all second countable completely 
regular spaces is equivalent to the category of subspaces 
of $\R^{\omega}$, which has only a set of objects\footnote{Such
a category is called a small category; however we do not want to
cause confusion with our later use of the word ``small.''}.
The only thing we have to remember is that such a category
is not closed under arbitrary products or coproducts.
But we will require that our category $\bT$ is closed
under finite limits and colimits (i.e. under finite products,
finite coproducts, pullbacks, and pushouts), and under taking
open subspaces.

Similar considerations apply to the category $\bM$ of smooth
manifolds, but we will assume, as usual, that
our manifolds are second countable Hausdorff spaces. Then
any manifold can be considered as a subspace of some
$\R^{N}$. So our category of smooth manifolds has only a set of objects.
Note that the category of smooth manifolds is
closed under finite products and coproducts, and taking open subspaces,
but not under pullbacks or pushouts.

None of these technical considerations will bother us in the sequel.
See Remark \ref{rem:whysmallcatOK} below.

We recall one basic categorical definition, familiar in many concrete 
situations. 
Given a category $\CC$ and an object $C$ of $\CC$, the \textit{category
of objects over $C$} (also called the slice category or  comma category) 
$\CC/C$ is the
category whose objects are arrows $f:D \ra C$ is $\CC$ and whose arrows
are commutative triangles
\[
  \xymatrix{ D' \ar[rr]^{g} \ar[dr]_{f'} &    & D \ar[dl]^{f} \\
                                         & C. &
  }
\]
There is an obvious forgetful functor $F:\CC/C \ra \CC$,
$F(D \stackrel{f}{\ra} C) = D$.

\begin{example}\label{ex:smalltopsite}
\textit{The small site of a topological space.} Let $X$
  be a topological space and let $\Op(X)$ be the
  category of open subsets of $X$, where arrows are
  given by inclusions of open sets. (Hence
  there is at most one arrow between any two objects.)
  Say that $\{U_{i} \ra U\}$ covers $U$ if $U = \bigcup U_{i}$
  (the usual definition of an open cover).
  This is easily seen to be a basis for a Grothendieck topology
  on $\Op(X)$. The covering sieve generated by 
  $\{U_{i}\}$ $U$ is the family of all sets $V$ such that
  $V \subset U_{i}$ for some $i$, i.e. the maximal refinement of
  $\{U_{i}\}$.

  The resulting site is called the
  \textit{small site} of the space $X$. This is the original
  and motivating example for the notion of a site.
  However it is special in that the underlying category
  is just a partial order; there are no nontrivial
  endomorphisms.
  
  If $X$ is a smooth manifold, we can treat it
  as a topological space and use the site $\Op(X)$.
\end{example}

\begin{example}\label{ex:largetopsite}
\textit{The large site of a topological space.} Let 
  $\bT$ be a fixed category of topological spaces
  as mentioned above, and let $X$ be a space in $\bT$.
  Let $\bTX$ be the category of spaces over $X$. 
  Define a basis on $\bTX$ by declaring 
  $\{f_{i}: Y_{i} \ra Y\}$ to be a covering family
  if each $f_{i}$ is an open embedding and
  $\bigcup_{i}f_{i}(Y_{i}) = Y$.
  Again this is easily seen to be a basis for a 
  Grothendieck topology on $\bTX$. The resulting
  site is called the \textit{large site} of the space
  $X$.
  
  An important special case is where we take $X = \ast$,
  a one-point space. Then $\bTX = \bT$. We will
  call this the \textit{large site of all topological spaces}
  or the \textit{absolute topological site} to distinguish
  it from the relative case of spaces over another space $X$.
\end{example}

Note that the small site of a space $X$ is the restriction of the
large site of $X$ to the full subcategory whose objects are
the open subsets of $X$. We will discuss the relationship of
these two sites (and their sheaves and stacks) at some length
in Sec.~\ref{sec:largesmall}.

\begin{example}
\label{ex:largemanifoldsite} 
  (The large site of a smooth manifold.) Let 
  $\bM$ be the category of smooth, second countable, Hausdorff manifolds
  as mentioned above, and let $M$ be a fixed smooth manifold in $\bM$.
  Let $\bMM$ be the category of manifolds over $M$.
  Define a basis on $\bMM$ by declaring 
  $\{f_{i}: N_{i} \ra N\}$ to be a covering family
  if each $f_{i}$ is an open embedding and
  $\bigcup_{i}f_{i}(N_{i}) = N$.
  This is a basis for a Grothendieck topology on $\bMM$. 
  The only nontrivial point is that the pullback
  of a diagram of manifolds
  \[
    \xymatrix{ & N_{1} \ar[d]^{f_{1}} \\ N_{2} \ar[r]_{f_{2}} & N_{3}}
  \]
  exists if one of the maps $f_{1}$ or $f_{2}$ is an open embedding.
  (In general, pullbacks do not exist in the category of smooth manifolds,
  but they do if one of the maps is a submersion, by transversality.)
  The resulting site is called the \textit{large site} of the manifold
  $M$. If necessary we will refer to it as the
  \textit{large smooth site} of $M$ to distinguish it
  from the topological site. Again we can take $M = \ast$,
  and we get \textit{the large site of all smooth manifolds}
  or the \textit{absolute smooth site}, which we will denote simply by $\bM$.
\end{example}
\begin{example}\label{ex:euclideansite}
  The site $\bM$ defined above is much larger than necessary for many
  purposes.
  If one wants to think of stacks as a replacement for the category
  of manifolds, it is more elegant not to have to define manifolds
  first. It is in fact possible to do this.

  Let $\mathbf{R}$ be the full subcategory of $\bM$  
  whose objects are $\{\R^{0}, \R^{1}, \R^{2}, \dots \}$,
  Define the \textit{Euclidean site} $\Euc \subset \bM$ 
  to be the full subcategory of $\bM$ obtained
  from $\mathbf{R}$ by taking open subsets and disjoint unions,
  equipped with the usual open cover topology (which is the same as 
  the restriction of the topology from Example \ref{ex:largemanifoldsite}. 
  We will see that this site is equivalent to the site $\bM$,
  in the sense of having the same sheaves and stacks 
  (Example~\ref{ex:manifoldcomparison}). In fact one could just use $\mathbf{R}$
  as the base site.

  Whether one takes $\mathbf{R}$, $\Euc$, or $\bM$ as the base site depends on
  one's goal. If one wants to redo manifold theory from the ground up,
  then taking $\R$ as the base site is most elegant, since it is the simplest
  site. However $\R$ is not the most convenient in practice, since
  it is not closed under the operations of disjoint union and taking open
  subsets. The site $\Euc$ is more convenient than $\R$, while still 
  not assuming anything
  about manifolds; however there are many natural constructions which
  jump out of $\Euc$, so the site $\bM$ is often most convenient. 
  Since we are happy to assume known results about manifolds, we 
  will usually use $\bM$ as the base site for maximal
  convenience. However we will occasionally take the foundational perspective,
  where we use $\Euc$ as the base site. See Example \ref{ex:manifoldassheaf} below.
\end{example}
\begin{example}\label{ex:alggeomsites}
  The original examples of sites come from algebraic geometry.
  One considers the category $\Sch$ of schemes and puts various
  topologies on it. The most basic topology is the \textit{Zariski
  topology}, where a covering family is a jointly surjective
  family of open embeddings. However the fundamental observation
  of Grothendieck which motivated the development of the theory
  of sites was that this topology is not fine enough
  to emulate the constructions one can do in the
  topological category. Rather, one should use a topology
  such as the \textit{\'etale topology}, where the maps in covering
  families need not be injective. We will not say much 
  about the algebraic geometry situation in this paper,
  as there is an extensive literature, e.g. \cite{SGA4.1}, 
  \cite{Giraud:Book}, \cite{LMB}, \cite{DeligneMumford}.
\end{example}

\begin{remark}\label{rem:whysmallcatOK}
  The above examples, especially the manifold example, give an idea
  why it is not a problem that we need to have our base category
  have a set of objects. The base category should be thought of
  as the category of local models (as in the case of $\Euc$).
  Stacks will essentially be spaces made by gluing these local
  models together, so they can be quite general even if the
  local models are simple.
\end{remark}

\subsection{Presheaves and Sheaves}\label{subsec:presheavessheaves}

We now turn to the standard definition of presheaves and sheaves over
a site. These generalize the familiar definitions for the case $\Op(X)$.
We will see in Secs.~\ref{subsec:fibgroupoid} and \ref{subsec:prestackstack}
that these are special cases of the definitions of prestacks and stacks,
for which we will use a slightly different language. We will
use the language of bases (covering families) as it is more familiar; 
see \cite{MacMoer:Sheaves} for a presentation in terms of covering sieves.

\begin{Def}\label{def:presheafsheaf}
  A \textbf{presheaf} $P$ on a category $\CC$ is a contravariant functor $P:\CC \ra \Sets$.
  We will often write the action of maps using restriction notation:
  for $f: C \ra D$ and $x \in F(V)$, we write
  \[
    F(U \stackrel{f}{\ra} V)(x) = x|_{U},
  \]
  the arrow $f$ being understood, when this will cause no confusion.
  If we need to be more careful we can use the functorial notation $F(f)(x)$.

  A \textbf{sheaf} $F$ on a site $(\CC,J)$ with basis $K$ is a presheaf on $\CC$
  such that, for every covering family $\{C_{i} \ra C\}$ in $K$ and every family of elements
  $x_{i} \in F(C_{i})$ such that $x_{i}|_{C_{ij}} = x_{j}|_{C_{ij}}$ for all $i,j$,
  there is a unique $x \in F(C)$ with $x|_{C_{i}} = x_{i}$.
\end{Def}
Here and in the rest of the paper we denote $C_{ij} = C_{i} \times_{C} C_{j}$.

Occasionally wants only the uniqueness part of the sheaf condition:
\begin{Def}\label{def:seppresheaf}
  A presheaf $P$ on a site $(\CC,J)$ with basis $K$ is \textbf{separated}
  if, given $x, x' \in F(C)$ and a covering family $\{C_{i} \ra C\}$ in $K$
  such that $x|_{C_{i}} = x'|_{C_{i}}$ for all $i$, we have $x = x'$.
\end{Def}
For example, on $\Op(X)$, the assignment 
$U \mapsto \{\text{bounded functions on~} U\}$ is a separated presheaf
which is not a sheaf.

We can present the data of a presheaf in another way, which will be the basis of the
generalization to stacks. Given a presheaf $F$ on $\CC$, we form the category
$\DD$ with set of objects
\[
  \Ob(\DD) = \coprod_{C \in \CC} F(C)
\]

and with arrows defined as follows. If $D \in F(C)$, $D' \in F(C')$,
and $f:C \ra C'$, then we create an arrow from $D$ to $D'$ if
$F(f)(D') = D$. Precisely, the set of arrows from $D$ to $D'$ is
\[
  \DD(D,D') = \{ f:C \ra C' \: | \: f(D') = D \}.
\]
To distinguish between the arrows of $\DD$ and those of $\CC$ we can use
a tilde: $\tilde{f}:D \ra D'$ is the arrow in $\DD$ corresponding to the arrow
$f:C \ra C'$ in $\CC$.
There is an obvious covariant functor from $\DD$ to $\CC$, taking
$D \in F(C)$ to $C$ and taking $\tilde{f}:D \ra D'$ to $f:C \ra C'$.

Clearly there is much that is special about this functor; we will see how
it is characterized in the next section.

\subsection{Fibered Categories}\label{subsec:fibgroupoid}

First we define a generalization of the notion of a presheaf.

\begin{Def}\label{def:fiberedgroupoid}
  Let $\CC$ be a category. A \textbf{category fibered in groupoids over $C$}
  is a category $\DD$ and a functor $F: \DD \ra \CC$
  such that:
  \begin{enumerate}
  \item \label{lift} Given any arrow $f: C' \ra C$ in $\CC$ and 
        an object $D$ of $\DD$ with $F(D) = C$, there is an
        arrow $g: D' \ra D$ with $F(g) = f$.
  \item \label{fillin} Suppose we have diagrams as follows (the top diagram
        is in $\DD$, the bottom in $\CC$):
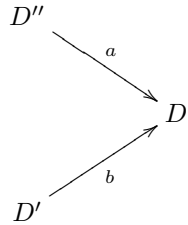
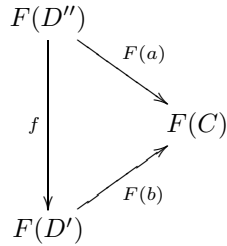
        \[
          \xymatrix{ D'' \ar[dr]^{a}        &       \\
                                            &  D    \\
                     D' \ar[ur]_{b}         &       \\
                                            &       \\
               F(D'') \ar[dr]^{F(a)} \ar[dd]_{f} &   \\
                                            &  F(C) \\
                     F(D') \ar[ur]_{F(b)}   &       \\
          }
        \]
        Then there is a unique $g: D'' \ra D'$ such that
        $b g = a$ and $F(g) = f$. (I.e. there is a unique
        way to fill in the top diagram such that its image
        under $F$ is the bottom diagram.)
  \end{enumerate}
\end{Def}
This is often referred to simply as a ``groupoid over $\CC$,''
but we will avoid that terminology to prevent confusion with
other notions of groupoid. We will however refer to it
simply as a ``fibered category over $\CC$,'' as we will
not have need of more general fibered categories.

Given an object $C$ of $\CC$, the \textit{fiber} $F_{C}$
is the subcategory of $\DD$ whose objects map to
$C$ and whose arrows map to $1_{C}$ under $F$.
It is easy to see that each fiber $F_{C}$ is indeed
a groupoid, i.e. every morphism in $F_{C}$ is invertible.

The most fundamental examples of fibered categories are the 
\textit{representables}, defined as follows. 
Let $C$ be an object of $\CC$ and consider the natural forgetful functor
$F: \CC/C \ra \CC$.
This makes $\CC/C$ into a fibered category over $\CC$. For,
given a map $f:D' \ra D$ and an object $g: E \ra C$ of $\CC/C$,
the composition $gf: D' \ra C$ is an object with $F(gf) = D'$,
and the map $f$ is an arrow from $gf$ to $g$.
In fact, this arrow is clearly the only arrow from $gf$ to $g$
whose image under $F$ is $f$. This makes condition
\ref{fillin} above easy to verify.
It also means that $\CC/C$ is fibered in
\textit{sets} over $\CC$ (we also say it is 
\textit{discretely fibered}): given any two objects
$g,g': D \stackrel{f}{\ra} C$ with the same source space
(i.e. with $F(g) = F(g')$, so $g,g' \in F_{D}$),
any arrow $a: D \ra D'$ with $F(a) = 1_{C}$ must be the identity.
Hence the fibers are sets, thought of as groupoids with no nonidentity
morphisms. We shall see shortly another, more familiar, way of describing
discretely fibered categories.

One can soup this up by replacing the base category
$\CC$ with $\CC/X$ for a fixed object $X$. Then any space
$h: C \ra X$ over $X$ defines a discretely fibered category 
$\CC/C \ra \CC/X$ by 
$F(D \stackrel{g}{\ra} C) = hg$.

\begin{example}\label{ex:representabletop}
  We can apply the foregoing to the category $\bT$ of spaces. 
  Let $Z$ be a space and consider the natural forgetful functor
  $F: \bTZ \ra \bT$. This is a fibered category, and the fibers are
  sets, as shown above. Or, we can apply the relative version,
  obtaining a fibered category $F: \bTZ \ra \bTX$ from
  any space over $X$, $h:Z \ra X$.
  
  Hence any space defines a discretely fibered category 
  over $\bT$; and any space over $X$ defines a discretely
  fibered category over $\bTX$.
  One can also do the parallel constructions in the smooth case.
  
  It is worth taking a moment to look at what this construction does.
  It describes a given space $Z$ (or a space $h:Z \ra X$ over a fixed base $X$)
  in terms of the maps into $Z$ from an arbitrary source space.
  We will see below that this is a slight rephrasing of a the standard
  construction of taking the functor which the space $Z$ represents.
\end{example}

\begin{example}\label{ex:principalbundles}
  Now we will look at the basic example of a fibered category
  with nondiscrete fibers. This is where genuine ``stackiness'' 
  enters. Let $G$ be a topological group (in $\bT$)
  and let $\Prin_{G}$ be the category whose objects are
  $G$-principal bundles $p:P \ra Y$, 
  and whose arrows from $p_{1}:P_{1} \ra Y_{1}$ to
  $p_{2}:P_{2} \ra Y_{2}$ are pairs 
  $(\tilde{f}: P_{1} \ra P_{2},f:Y_{1} \ra Y_{2})$,
  where $\tilde{f}$ is a $G$-equivariant map and 
  $f p_{1} = p_{2} \tilde{f}$. (This is the same
  data as a map $f$ on the base spaces and a map $P_{1} \ra f^{*}P_{2}$.)
  There is an obvious forgetful functor $F:\Prin_{G} \ra \bT$.
  It is easy to see that the pullback operation on principal
  bundles makes $\Prin_{G}$ into a fibered category.
  The fiber over a space $Y$ is the usual category of
  $G$-principal bundles over $Y$; it is a well-known
  fact that any morphism of principal bundles over $Y$
  is an isomorphism, so this is indeed a groupoid.

  Again the smooth case is exactly parallel.
  We will see below (Section~\ref{subsec:assocstack})
  that this is a special case of taking the stack quotient
  of a group action.
\end{example}

We can think of a fibered category $F: \CC \ra \DD$ as defining
a ``presheaf of groupoids'' $\F$ on $\CC$ in the following way. 
First we need to choose, for every
arrow $f: C' \ra C$ in $\CC$ and for every object $D$ in $F_{C}$,  
a lift $\tilde{f}:f^{*}D \ra D$ of $f$ as guaranteed by axiom \ref{lift}
of Definition~\ref{def:fiberedgroupoid}. We will speak
of $\tilde{f}$ as a \textit{pullback arrow} for $f$ and
call $f^{*}D$ the \textit{pullback} of $D$ by $f$.
When it will not cause confusion we will also refer to
$f^{*}D$ as $D|_{C'}$.

Now we define $\F$. 
To every object $C$ on $\CC$ we assign the fiber $\F(C) = F_{C}$, which
is a groupoid.  To every arrow $f: C' \ra C$ we need to assign a functor
$\F(f)$ from $F_{C}$ to $F_{C'}$. On objects, this functor is defined
by $\F(f)(D) = f^{*}D$. 
Given an arrow $g: D_{1} \ra D_{2}$ in $F_{C}$, 
$\F(f)(g)$ is the unique arrow 
$f^{*}g: f^{*}D_{1} \ra f^{*}D_{2}$ (guaranteed by axiom \ref{fillin}) 
which fits into the diagram
\[
  \xymatrix{ f^{*}D_{1} \ar[r] \ar[d]_{f^{*}g} & D_{1} \ar[d]^{g} \\
             f^{*}D_{2} \ar[r]                 & D_{2} 
  }
\]
where the horizontal arrows are pullback arrows.
Again, where appropriate we will refer to $f^{*}g$ as $g|_{C'}$.

This does not quite define a presheaf of groupoids, for the
following reason: we generally have $\F(f)\F(f') \ne \F(f' f)$.
So $\F$ is not strictly functorial.
However it is easy to show, using axiom \ref{fillin},
that the pullback functors $\F(f)$ satisfy 
\[
  \F(f) \circ \F(f') \stackrel{\phi_{f,f'}}{\iso} \F(f' f)
\]
and that the maps $\phi_{f,f'}$ satisfy an appropriate coherence condition.
(Roughly, any diagram involving $\phi$ and $\F$ which ought to commute,
does.) This is described by saying 
that $\F$ defines a \textit{lax} functor from
$\CC$ to the category of groupoids, i.e. a lax presheaf of groupoids.

Given a different choice of pullback arrows, one gets a 
different, but equivalent, lax presheaf of groupoids. 
(The appropriate notion of equivalence is left to the reader.)
One can also go the other way, from a given lax presheaf of groupoids
to a fibered category. 
We will use both points of view on fibered categories, both
the definition and the characterization as lax
presheaves of groupoids. The language of lax presheaves of
groupoids is often more intuitive (by analogy with ordinary
presheaves) but explicitly dealing with the isomorphisms
$\phi_{f,f'}$ can be annoying. The language of fibered categories
is more elegant and avoids both choosing pullback arrows and
dealing with the $\phi_{f,f'}$. See \cite{Brylinski:LoopBook}
for a treatment using the language of (lax) presheaves of groupoids.

Regardless of the language chosen,
the laxness does not usually present a problem, except in stating
things carefully. However, sometimes we will use
\textit{strict} presheaves of groupoids (which are functorial on the nose).
The corresponding fibered category will also be called strict.
(See Proposition~\ref{prop:strictify} below in this context.)

\begin{example}\label{presheaves}
\textit{Ordinary Presheaves}. Let $F: \CC \ra \DD$ be a discretely
fibered category. Then it is easy to see that pullback
arrows (and hence objects) are uniquely defined:
given $f: C' \ra C$ in $\CC$ and $D$ an object of $F_{C}$,
there is a unique arrow $\tilde{f}: f^{*}D \ra D$ in $\DD$ 
with $F(\tilde{f}) = f$. The uniqueness implies that
$\F(f) \circ \F(f') = \F(f' f)$, so $\F$ defines an
honest contravariant functor from $\CC$ to $\Sets$,
i.e. a presheaf of sets. 

It is easy to see that this inverts the construction from Sec.~\ref{subsec:presheavessheaves},
showing that presheaves are equivalent to discretely fibered categories.
We will use this equivalence throughout the paper.
\end{example}

\begin{example}\label{ex:reppresheaves}
\textit{Representable Presheaves}
  Let $Z$ be a space and consider the natural forgetful functor
  $F: \bTZ \ra \bT$ as a discretely fibered category, as we did above.
  The fiber above each space $X$ is the set $\bT(X,Z)$ of
  continuous maps from $X$ to $Z$. The pullback functor is easily
  seen to be given by composition of maps. Hence the presheaf
  defined by $F$ is just the presheaf $\uZ = \Maps(-,Z)$ 
  represented by 
  $Z$. The Yoneda Lemma asserts that any space $Z$
  is determined (up to canonical isomorphism) by the presheaf $\uZ$.
  Hence we can think of a fibered category as a very loose
  kind of ``generalized space.'' 
  (This will be further justified in Sec.~\ref{sec:stackmaps} below.)
  Further, it is easy to verify that $\uZ$ is in fact a sheaf.
  (This is just the gluing lemma for continuous maps.)
  We will come back to this point shortly.
\end{example}

More generally, in any category $\CC$, the
representable discretely
fibered category defined by an object $C$ 
is equivalent to the usual representable functor
(presheaf of sets) $\underline{C} = \Hom_{\CC}(-,C)$.

\begin{example}\label{ex:principalbundleslax}
  \textit{Principal Bundles.} Recall the fibered category 
  $\Prin_{G}$ of principal bundles over the category of spaces $\bT$.
  The usual pullback operation for bundles by a map $f: Y \ra X$, given by
  \[
    f^{*}P = Y \times_{X} P
  \]
  defines a pullback object, and the universal arrow 
  $\tilde{f} = \pi_{2}: Y \times_{X} P \ra P$ defines a pullback arrow.
  (This is in fact one of the motivations for using the term ``pullback.'')
  For a composition $f' f$ of maps of spaces, the corresponding
  pullbacks are canonically isomorphic but not identical, hence this
  defines a lax functor, not a strict one.
\end{example}

One use of the presheaf language is the following.
Suppose we have a fibered category $F: \DD \ra \CC$, 
with a choice of pullback arrows, an object $C$ of $\CC$, 
and two objects
$D_{1}, D_{2}$ of $F_{C}$. Then we 
define the \textit{presheaf of local isomorphisms}
$\IIso(D_{1},D_{2})$ from $D_{1}$ to $D_{2}$ as follows. 
Given an arrow $f:C' \ra C$ in $\CC$, let
\[
  \IIso(D_{1},D_{2})(C' \stackrel{f}{\ra} C) 
     = \Iso_{F_{C'}}(f^{*}D_{1},f^{*}D_{2}).
\]
It is easy to check (using axiom \ref{fillin}) 
that this is a presheaf (of sets).

\subsection{Prestacks and Stacks}\label{subsec:prestackstack}

Now we will give one of the conditions for a fibered category to
be a stack.
\begin{Def}\label{def:prestack}
  Let $(\CC,J)$ be a site and 
  let $F:\DD \ra \CC$ be a fibered category. $F$ is a
  \textbf{prestack} if for every object $C$ of $\CC$
  and every two objects $D_{1},D_{2}$ of $F_{C}$,
  the presheaf $\IIso(D_{1},D_{2})$ is a sheaf over the
  category $\CC/C$.
\end{Def}
(The Grothendieck topology on $\CC/C$ is the one induced
by $J$; we leave the easy definition to the reader.)

\begin{example}\label{ex:principalprestack}
  The fibered category $\Prin_{G}$ of principal
  bundles is easily seen to be a prestack; this is just the gluing axiom for 
  maps of spaces, applied to principal bundles. 
\end{example}

\begin{example}\label{ex:separatedpresheaf}
  Let $F: \DD \ra \CC$ be a discretely fibered category, 
  i.e. a presheaf of sets. In this case two objects are isomorphic if and only
  if they are equal. Suppose we are given two objects 
  $D_{1},D_{2} \in F_{C}$ and a cover $\{ C_{\alpha} \ra C\}$.
  We will have consistent local sections of $\IIso(D_{1},D_{2})(C_{\alpha})$
  exactly when $D_{1}|_{C_{\alpha}} = D_{2}|_{C_{\alpha}}$ for every $\alpha$.
  Hence $F$ is a prestack if and only if the corresponding presheaf is 
  \textit{separated}. 
  Hence we can think of prestacks as a generalization of separated
  presheaves. 
\end{example}

The foregoing definition is the most commonly seen one;
however we can phrase it without resorting to a fixed choice of
pullbacks, as in the following.

Let $(\CC,J)$ be a site and let $F: \DD \ra \CC$ be a fibered category.
Let $C \in \CC$, let $\{C_{\alpha} \ra C\}$ be a cover of $C$ and let $x,y \in F_{C}$.
Suppose we have a diagram
\begin{equation}\label{firstprediag}
\xymatrix@=5pt{
     & & x_{\beta} \ar[dr] \ar[ddd] & 
      & & & & & & C_{\beta} \ar[dr] \ar@{=}[ddd] & \\
  x_{\alpha \beta} \ar[urr] \ar[dr] \ar[ddd] & & & x \ar@{.>}[ddd]^{\psi} 
    & & & & C_{\alpha \beta} \ar[urr] \ar[dr] \ar@{=}[ddd] & & & C \ar@{=}[ddd] \\
     & x_{\alpha} \ar[urr] \ar[ddd] & & & \ar@{~>}[rr]^{F}
       & & & & C_{\alpha} \ar[urr] \ar@{=}[ddd] & & \\
     & & y_{\beta} \ar[dr] & 
      & & & & & & C_{\beta} \ar[dr] & \\
  y_{\alpha \beta} \ar[urr] \ar[dr] & & & y 
   & & & & C_{\alpha \beta} \ar[urr] \ar[dr] & & & C \\
     & y_{\alpha} \ar[urr] & & 
     & & & & & C_{\alpha} \ar[urr] & &
  }
\end{equation}
We note some conventions:
we indicate by the squiggly arrow that
the diagram on the left gets mapped to the diagram on the right
under $G$; and the dotted arrow, as is usual, signifies
a desired fill-in arrow which we do not yet have.
We will avoid naming arrows in these kinds of diagrams when it is not necessary.
As usual we abbreviate $C_{\alpha} \times_{C} C_{\beta}$ by $C_{\alpha \beta}$, etc.
In this diagram the objects $\{x_{\alpha}\}$ and $\{x_{\alpha \beta}\}$ are arbitrary;
they play the role of the pullbacks of $x$ to $C_{\alpha}$, $C_{\alpha \beta}$ etc.

\begin{lemma}\label{lem:prestackwithoutpullbacks}
  The fibered category $F: \DD \ra \CC$ 
  is a prestack if and only if, given a diagram of the form (\ref{firstprediag}),
  there is a unique arrow $\psi :x \ra y$ filling in the dotted arrow
  in the left hand diagram in (\ref{firstprediag}), with $F(\psi) = 1$.
\end{lemma}
The proof of this lemma is straightforward and we will omit it. The point
is that \textit{any} arrow $x_{\alpha} \ra x$ serves as a pullback,
when $F: \DD \ra \CC$ is a category fibered in groupoids, and the consistency
condition on overlaps is expressed by asserting that there is \textit{some}
object $x_{\alpha \beta}$ over $C_{\alpha \beta}$ over which the maps agree.
We will use this characterization of prestacks in Sec.~\ref{sec:stackoverstack}.

Before we give the definition of a stack we need to define descent data. 
\begin{Def}\label{def:descentdatawithpullbacks}
  Let $F: \DD \ra \CC$ be a fibered category,
  with a choice of pullback arrows, and
  let $C$ be an object of $\CC$ and let $\mathcal{C}: f_{i}: C_{i} \ra C$, 
  be a covering family of $C$. Let $C_{ij} = C_{i} \times_{C} C_{j}$
  and $C_{ijk} = C_{i} \times_{C} C_{j} \times_{C} C_{k}$.
  Then \textit{descent data} for $F$
  over $\mathcal{C}$ are a family of objects $D_{i}$ of $F_{C_{i}}$
  and arrows 
  $\phi_{ji}: D_{i}|_{C_{ij}} \stackrel{\sim}{\ra} D_{j}|_{C_{ij}}$ 
  in $F_{C_{ij}}$
  (necessarily isomorphisms) such that 
  \begin{equation}\label{cocycle}
    (\phi_{ij}|_{C_{ijk}}) \circ (\phi_{jk}|_{C_{ijk}}) 
       = \phi_{ik}|_{C_{ijk}}.
  \end{equation}
  We will always use descent data that is \textit{normalized}, meaning
  that $\phi_{ii} = 1$.
\end{Def}

\begin{example}\label{ex:principalstack}
  For $\Prin_{G}$, descent data over a space
  $X$ and an open cover $\UU = \{ U_{i} \subset X \}$ is just the
  usual data of a principal bundle $P_{i}$ over each $U_{i}$ and transition
  functions on the overlaps which satisfy the usual cocycle
  condition (\ref{cocycle}). Given such descent data, one knows
  that there exists a principal bundle $P$ over $X$ which
  is isomorphic to $P_{i}$ over each $U_{i}$, in a way that is
  compatible with the transition maps. In the terminology
  of Grothendieck, we say that all such descent data are \textit{effective}.
  This is the property that we use as the definition of a stack.
\end{example}

\begin{Def}\label{def:stack}
  Let $(\CC,J)$ be a site and let $F: \DD \ra \CC$ be a prestack. 
  Then $F$ is a \textbf{stack}
  if all descent data $(D_{i},\phi_{ij})$ for $F$ 
  with respect to a covering family $f_{i}: C_{i} \ra C$ are effective,
  that is, given such descent data, there is an object $D$ of $F_{C}$
  and isomorphisms $\psi_{i}: D|_{C_{i}} \rai D_{i}$ such that 
  \[
    \psi_{j}|_{C_{ij}} = \phi_{ji} \circ \psi_{i}|_{C_{ij}}. 
  \]
\end{Def}
Since $F$ is a prestack, this object $D$ is unique up to 
a canonical isomorphism. 

\begin{example}\label{ex:sheafasstack}
  \textit{Sheaves as stacks.} It is easy to verify that a presheaf (as
  a discretely fibered category) is a stack if and only if it is
  a sheaf. In particular, let $Z$ be a space in $\bT$
  and consider the fibered category $F: \bTZ \ra \bT$,
  or in other words the representable functor $\uZ$.
  As mentioned above this is a sheaf, so it defines a stack.
  This is the sense in which we will think of stacks as
  generalized spaces.
\end{example}

More generally, let $X$ be a space and let $h:Z \ra X$ be
a space over $X$. Consider the discretely fibered category
$F: \bTZ \ra \bTX$ as above, or in other words the
representable functor $\uZ$ defined by $Z$ (as a space over $X$).
Once again this is a sheaf, hence a stack. So we will think
of stacks over the large site $\bTX$ as generalized spaces
over $X$.

Even more generally, consider a general site $(\CC,J)$.
Suppose that $J$ is \textit{subcanonical}, i.e. that
every representable presheaf is a sheaf. (This is almost
always true in examples, in particular, it will be true
for all examples of sites in this paper. So we will assume it henceforth
when necessary.) Then we can consider the
presheaf $\uC$ represented by any object $C$ as a stack.
We will use this constantly throughout the paper.

\begin{example}\label{ex:smallsheaf}
  \textit{Small sheaves on a space.}
  On the other hand we can consider the small site
  $\Op(X)$. A sheaf over this site is a sheaf of sets in the ordinary
  sense. This can be considered as a space over $X$ by constructing
  the \'etale space of the sheaf. However the spaces we get
  are clearly a proper subcategory of the category of all
  spaces over $X$, much less all sheaves over the large site $\bTX$.
  So we need to be clear about which site we mean when we refer
  to a ``sheaf over X.'' We will discuss this distinction,
  and the corresponding one for stacks, in Section~\ref{sec:largesmall}.
\end{example}

\begin{example}\label{ex:prinstackagain}
  We saw above that $\Prin_{G}$ is a stack,
  because we can glue principal bundles with compatible transition
  maps. We will see below that it is a somewhat special stack
  (a gerbe), but it is still the best example to keep in mind
  when thinking about stacks.
\end{example}

Note that a fibered category could be thought of as a
``pre-prestack,'' that is, a stack with none of the sheaflike
conditions included. We will make some definitions for
arbitrary fibered categories, since they are not harder in this
case. The relationship of fibered categories to stacks is analogous
to that between presheaves and sheaves. In particular there is
a canonical way to associate a stack to a fibered category.
However it turns out that many constructions automatically
yield prestacks, and one only has to ``stackify'' these, instead
of starting with a fibered category that is not a prestack.
We will explain this construction in the next section.

As with the prestack condition, a fixed choice of pullbacks is
not necessary to express the stack condition, though it
is usually written that way. Analogous to 
Lemma \ref{lem:prestackwithoutpullbacks},
we have the following.

Let $(\CC,J)$ be a site and let $F: \DD \ra \CC$ be a prestack.
Let $C \in \CC$ and let $\{C_{\alpha} \ra C\}$ be a cover of $C$.
Descent data in this context, without a fixed choice of pullbacks,
is expressed as follows.
For each $\alpha$ let $x_{\alpha} \in G_{C_{\alpha}}$,
for every pair $\alpha, \beta$ let there be a diagram
\begin{equation}\label{zerothstackdiag}
  \xymatrix@=5pt{     
      & x_{\alpha} & & & & & C_{\alpha} \\
    x_{\alpha \beta} \ar[ur] \ar[dr] & & \ar@{~>}[rr]^{F} & & 
        & C_{\alpha \beta} \ar[ur] \ar[dr] & \\
      & x_{\beta}  & & & & & C_{\beta} 
  }
\end{equation}
(where the arrows on the right are the canonical projections)
and for every triple $\alpha ,\beta ,\gamma$ let there be a commutative diagram
\begin{equation}\label{firststackdiag}
  \xymatrix@=5pt{ 
    & x_{\alpha \beta} \ar[r] \ar[dr] & x_{\alpha} 
      & & & & & C_{\alpha \beta} \ar[r] \ar[dr] & C_{\alpha} \\
    x_{\alpha \beta \gamma} \ar[ur] \ar[r] \ar[dr] 
       & x_{\alpha \gamma} \ar[ur] \ar[dr] & x_{\beta}
      & \ar@{~>}[rr]^{F} & & & C_{\alpha \beta \gamma} \ar[ur] \ar[r] \ar[dr] 
       & C_{\alpha \gamma} \ar[ur] \ar[dr] & C_{\beta} \\
    & x_{\beta \gamma} \ar[ur] \ar[r] & x_{\gamma}
      & & & & & C_{\beta \gamma} \ar[ur] \ar[r] & C_{\gamma} 
  }
\end{equation}
where again all of the maps on the right are the canonical projections. 

\begin{lemma}\label{lem:stackwithoutpullbacks}
  Let $F: \DD \ra \CC$ be a prestack. 
  $F$ is a stack if and only if the following condition is satisfied
  for every cover $\{C_{\alpha} \ra C\}$ in $\CC$.
  Given $x_{\alpha} \in G_{C_{\alpha}}$ and
  diagrams of the form (\ref{zerothstackdiag})
  and (\ref{firststackdiag}), there is an $x \in G_{C}$
  and arrows $\{x_{\alpha} \ra x\}$ filling in the commutative diagram
  \begin{equation}\label{firststackwantdiag}
    \xymatrix@=5pt{  
        & x_{\alpha} \ar[dr] & & & & & & C_{\alpha} \ar[dr] &  \\
      x_{\alpha \beta} \ar[ur] \ar[dr] & & x 
       & \ar@{~>}[rr]^{F} & & & C_{\alpha \beta} \ar[ur] \ar[dr] & & C \\
        & x_{\beta} \ar[ur]  & & & & & & C_{\beta} \ar[ur] &  
    }
  \end{equation}
\end{lemma}

To conclude this section we include one more characterization
of descent data in terms of covering sieves. Given a fibered category
$F: \DD \ra \CC$, an object $C$ of $\CC$, let $S$ be a covering sieve of $C$.
Note that $S$ can be considered as a full subcategory of $\CC/C$.
Denote by $i: S \ra \CC$ the forgetful functor. 
Then \textit{descent data for $F$ over the sieve $S$} is
a functor $A: S \ra \DD$ such that $F \circ A = i$.
In other words, we assign an object of $F_{C'}$ to every arrow
$C' \ra C$ in $S$, and a family of compatible pullback arrows connecting
all of these objects. In contrast to the previous definition of
descent data (Definition~\ref{def:descentdatawithpullbacks}),
where we used a minimal amount of data (objects defined only for
a covering family), this new definition uses the maximal amount of
data (objects defined for the whole covering sieve). So it
is less computationally useful but more elegant, requiring no
choices. It is not hard to show that this notion of
descent data is equivalent to our previous definitions.

Note that descent data for $F$ over a given sieve $S$ form a category
(in fact a groupoid) $\Desc(S,F)$, 
the morphisms being natural transformations $\alpha: A \ra B$
(which are necessarily isomorphisms).
Also, if we have an inclusion of sieves $S \subset T$,
there is a natural restriction functor $\Desc(T,F) \ra \Desc(S,F)$.

One can use this notion of descent data to give an
elegant version of the definition of prestacks and stacks;
see \cite{Giraud:Book}. We will use this version of descent
data in defining the stack associated to a prestack
(Def. \ref{def:stackification}).

\subsection{Maps of Stacks}\label{sec:stackmaps}

One advantage of thinking of stacks as fibered categories (as opposed to
lax sheaves of groupoids) is 
that the notion of a map of stacks is very simple in this language.
However we will immediately see that 2-categorical notions
(see e.g. \cite{CWM}) appear, and they are essential to understanding stacks fully.

\begin{Def}\label{def:mapoffiberedcats}
  Let $F: \DD \ra \CC$ and $F': \DD' \ra \CC$ be 
  two fibered categories over $\CC$. Then a \textbf{map of fibered categories}
  (or simply, a map) $a$ from $F'$ to $F$ is a functor $a: \DD' \ra \DD$
  such that $F \circ a = F'$. We will sometimes call a map a
  \textbf{1-arrow} to distinguish it from the following.

  A \textbf{2-isomorphism} from a map $a: F' \ra F$ to a map $a':F' \ra F$
  is a natural isomorphism of functors.
\end{Def}
We will denote the groupoid (maps, 2-isomorphism between maps) between
two fibered categories $F', F$ by $\HOM(F',F)$. The set of maps alone will
be denoted by $\Hom(F',F)$. Given any map of fibered
categories $a: F' \ra F$ and any object $C$ of the base
category $\CC$, we get a functor $a_{C}:F'_{C} \ra F_{C}$
by restriction.

A map of stacks (or prestacks) is just a map of the underlying 
fibered categories, and similarly for 2-isomorphism.
Denote by $\St(\CC,J)$ (or $\St(\CC)$ by abuse of notation)
the 2-category of stacks over the site $(\CC,J)$.

We leave it to the reader to characterize maps of fibered categories
in the language of lax presheaves of groupoids. However it is important
to note that a map of fibered categories will not in general
take pullbacks to pullbacks (on the nose), as these are extra data.
Instead one will get preservation of pullbacks only up to
coherent isomorphism. If we stick to the fibered category language
we can avoid worrying about these isomorphisms and their particular
coherence conditions.

Just as one rarely has two categories being isomorphic,
we will not often see a strict isomorphism of stacks. Rather,
we say an \textit{equivalence} of stacks is a map $\phi:F' \ra F$
with a quasi-inverse $\psi:F \ra F'$, i.e. we have 2-isomorphisms
$\phi \psi \iso 1$ and $\psi \phi \iso 1$.
Just as one usually treats equivalent categories
as ``the same,'' we will consider equivalent stacks
as containing the same information. 

The following lemma follows easily from the definition
of a fibered category.
\begin{lemma}\label{lem:equiv}
  Let $F: \DD \ra \CC$ and $F':\DD' \ra \CC$ be two fibered 
  categories and let $a:F' \ra F$ be a functor with $Fa = F'$.
  Then $a$ is an equivalence of fibered categories (i.e.
  it has a quasi-inverse which also respects $F,F'$)
  if and only if for each object $C$ of $\CC$, 
  the functor $a_{C}:F'_{C} \ra F_{C}$ is an equivalence of 
  categories.
\end{lemma}

We will also use the following easy lemma.
\begin{lemma}\label{lem:equivtostackisstack}
  Let $F,G$ be two fibered categories and let $\phi: F \ra G$
  be an equivalence of fibered categories. If $F$ is a stack (resp. prestack)
  then $G$ is a stack (resp. prestack).
\end{lemma}

\begin{example}\label{ex:stacknoauts}
  Let $F: \DD \ra \CC$ be a stack such that
  for every object $C$ of $\CC$ and for every object
  $D$ of $F_{C}$, the group of automorphisms 
  $\Iso(D,D)$ consists only of the
  identity. Let $\EE$ be the category formed from
  $\DD$ be identifying any two objects in each fiber which
  are isomorphic. (This kind of quotient of a category does
  not work in general, but it makes sense
  in this case because of the lack of nontrivial automorphisms.)
  There is a natural functor $G:\EE \ra \CC$ which is
  a discretely fibered category, and we have a 
  natural map of fibered categories $a:F \ra G$ which is
  an equivalence on every fiber. By Lemma~\ref{lem:equiv},
  $a$ is an equivalence of fibered categories.
  Hence any stack with only trivial automorphism groups $\Iso(D,D)$
  is equivalent to a sheaf.
\end{example}

We want to think of stacks over $\CC$ as generalized objects of $\CC$,
and this is made explicit by the notion of the representable stack $\uC$
associated to an object $C$ of $\CC$, as in Example \ref{ex:reppresheaves}.
So we need to check that the notions of map are compatible.
The following lemma is a straightforward consequence of the Yoneda Lemma.

\begin{lemma}\label{lem:mapsofreps}
  Let $C$ and $C'$ be objects of a subcanonical site $(\CC,J)$ and
  let $\uC$ and $\uC'$ be the corresponding representable stacks.
  Given an arrow $a: C' \ra C$, we get a map 
  $\underline{a}: \uC' \ra \uC$
  of the stacks by 
  $\underline{a}(B \stackrel{g}{\ra} C') = (B \stackrel{ag}{\ra} C)$.
  Further, every map $b: \uC' \ra \uC$
  of stacks comes from an arrow from $C'$ to $C$ in this way.
  Also, since $\uC$ and $\uC'$ are discretely fibered, there are no nontrivial
  2-isomorphisms between maps from $\uC'$ to $\uC$.
  Put more concisely, we get a full embedding of 2-categories 
  $\mathbf{y}: \CC \ra \St(\CC,J)$.
\end{lemma}

\begin{example}\label{ex:mapsofrepresentables}
Let $Z$ and $Z'$ be spaces in $\bT$ and consider
the corresponding representable stacks $\uZ$ and $\uZ'$. 
Given a continuous map $a: Z' \ra Z$, we get a map 
$\underline{a}: \uZ' \ra \uZ$
of the stacks by 
$\underline{a}(Y \stackrel{g}{\ra} Z') = (Y \stackrel{ag}{\ra} Z)$,
and every map $b: \uZ' \ra \uZ$
of stacks comes from a continuous map from $Z'$ to $Z$ in this way.
So we get a full embedding of 2-categories 
$\mathbf{y}: \bT \ra (\text{Topological Stacks})$, justifying
our viewpoint that topological stacks are generalized spaces.
\end{example}

In fact a map from a representable stack into an arbitrary
stack is always simple. The following lemma is a simple extension 
of the Yoneda Lemma.
\begin{lemma}\label{lem:mapfromrepresentable}
  Let $(\CC,J)$ be a subcanonical site. Let $C$ be an object of $\CC$
  and let $\uC$ be the stack represented by $C$.
  and let $F$ be a stack over $\CC$.
  Then there is a canonical equivalence of categories (groupoids)
  \[
    F_{C} \rai \HOM(\uC,F).
  \]
\end{lemma}
This means we should think of the fiber $F_{C}$ of $F$ over
$C$ as consisting of the maps from $C$ to $F$, where $F$ is considered
as a generalized space. In fact we will usually write 
$\HOM(C,F)$ instead of $\HOM(\uC,F)$. But we must pay attention to the
fact that $\HOM(\uC,F)$ is a groupoid, not just a set.

Lemma~\ref{lem:mapfromrepresentable} is crucial for intuition about
stacks. It continues the idea that we can understand a space by looking
at all maps into that space. 

\begin{example}\label{ex:manifoldassheaf}
We take the foundational perspective here and show how one can
think of a smooth manifold as a sheaf on the site $\Euc$ of ``local models.''
Let $M$ be a smooth manifold and consider the sheaf
$\uM$ represented by $M$, on the site $\bM$. This describes $M$ in
terms of all possible maps from all possible manifolds into $M$.
This is highly redundant, and rather tautological; 
we can cut this information down as follows. 
Consider the Euclidean site $\Euc$ defined in Example~\ref{ex:euclideansite},
which is a subcategory of $\bM$. We will use $U$ to denote
an object of $\Euc$, i.e. a disjoint sum $U = \coprod_{\alpha} U_{\alpha}$, where
$U_{\alpha} \subset \R^{n_{\alpha}}$.

The restriction of $\uM$ to $\Euc$, i.e. 
$\uM: U \mapsto C^{\infty}(R^{n},M)$ is a sheaf on $\Euc$.
Further any map $f:M \ra N$ gives a map $\uf:\uM \ra \uN$
on $\bM$, hence on $\Euc$. Note that $\uM(\R^{0})$ is the
underlying set of $M$, and $\uf$ restricted to $R^{0} \in \Euc$
is exactly the map $f$ as a map of sets. 
Since any map of manifolds is determined
by the corresponding point-set map, the map $\uf$ of sheaves
on $\Euc$ determines $f$. Since $M$ is covered by coordinate charts,
it is easy to show directly that
any map of sheaves $\phi: \uM \ra \uN$ comes from a smooth map
of manifolds $f:M \ra N$.

Hence the category of smooth manifolds is a full subcategory of the
category of sheaves on $\Euc$. 
\end{example}

In fact this is an example of the \textit{comparison lemma}, which
we present below. First we give a definition which will be useful
elsewhere as well, of the \textit{direct image} of a stack.

\begin{DefProp}[\cite{Giraud:Book}]\label{def:directimage}
  Let $F: \DD \ra \CC$ be a fibered category and let $u: \bA \ra \CC$
  be a functor. The \textbf{direct image} of $F$ under $u$, denoted
  $u_{*}F$,
  is the fibered category given by the canonical projection 
  \begin{equation}\label{directimage}
    \DD \times_{\CC} \bA \ra \bA.
  \end{equation}
  Further, suppose $J$ is a topology on $\CC$ and $\bA$ is given the
  induced topology. 
  If $F$ is a stack (resp. prestack) then $u_{*}F$ 
  is a stack (resp. prestack).
  If $F$ is a sheaf (resp. presheaf) then $u_{*}F$ 
  is a sheaf (resp. presheaf),
  which agrees with the usual direct image of sheaves (resp. presheaves).
\end{DefProp}
\begin{note}
  In (\ref{directimage}), the pullback of categories has
  $\Ob(\DD \times_{\CC} \bA) = \Ob(\DD) \times_{\Ob(\CC)} \Ob(\bA)$,
  $\Ar(\DD \times_{\CC} \bA) = \Ar(\DD) \times_{\Ar(\CC)} \Ar(\bA)$.
  So in particular, the fibers satisfy
  \[
    (u_{*}F)_{A} = F_{u(A)}.
  \]
  If we think of $F$ as a lax presheaf $\F$ of groupoids, then 
  $u_{*}F$ corresponds to the composed presheaf $\F \circ u$.
  (Note that this proves the last statement of the proposition.)
\end{note}

\begin{remark}\label{rem:directnotinverse}
  One may wonder why, since the direct image takes a stack
  over $\CC$ to a stack over $\bA$, it is not called the ``inverse image.''
  This is because the functor $u: \bA \ra \CC$ is regarded as a
  map of sites from $(\CC,J)$ to $(\bA,J_{A})$. This in turn comes from
  the case of topological spaces, where a continuous map $f: X \ra Y$
  induces a functor $f^{-1}: \Op(Y) \ra \Op(X)$. So maps of sites
  have the same ``direction'' as maps of topological spaces. 
\end{remark}

\begin{example}\label{ex:restrictedstack}
  Let $(\CC,J)$ be a site and let $i: \bA \subset \CC$ be a full subcategory.
  Then the direct image $i_{*}F$ of a stack on $\CC$ is also
  denoted $F|_{\bA}$ and is called the \textit{restriction} of $F$
  to $\bA$.   
\end{example}

Now we present the comparison lemma.
\begin{lemma}\label{lem:comparison}
  Let $(\CC,J)$ be a site and let $\bA \subset \CC$ be a full subcategory.
  Assume that the inclusion functor $i: \bA \subset \CC$ preserves all 
  pullbacks that exist in $\bA$. Also, suppose that
  every object $C$ of $\CC$ can be covered by objects from $\bA$.
  Then restriction to $\bA$ $F \ra i_{*}F$ induces an equivalence of
  categories between $\Sh(\CC)$ and $\Sh(\bA)$, and
  an equivalence of $2$-categories between $\St(\CC)$ and $\St(\bA)$.
\end{lemma}
\begin{note}
  The equivalence of $2$-categories in the lemma is a weak notion.
  In particular every stack $F$ on $\CC$ is \textit{equivalent}
  to a stack coming from a stack on $\bA$, but not necessarily isomorphic
  to such a stack. This will not bother us.
\end{note}
\begin{proof}
  The comparison lemma for sheaves is standard, see \cite{SGA4.1}, p. 288,
  or \cite{MacMoer:Sheaves}, p. 588. One also has the result of 
  Giraud \cite{Giraud:Book}, p. 91, which says that the $2$-category
  of stacks on a site depends (up to equivalence, as noted above)
  only on the category of sheaves over the site. Hence the
  comparison lemma for sheaves implies the one for stacks.  
\end{proof}

\begin{example}\label{ex:ordinarybasis}
  The most intuitive example of the use of the comparison lemma
  is in the small site $\Op(X)$ of a topological space $X$.
  If $\BB \subset \Op(X)$ is a basis (in the ordinary sense)
  for the topology of $X$, then the comparison lemma applies,
  and any small sheaf or stack on $X$ is determined by its
  values on $\BB$, as one would expect.
\end{example}

\begin{example}\label{ex:manifoldcomparison}
  We return to manifolds, continuing Example~\ref{ex:manifoldassheaf}.
  The inclusion $i:\Euc \ra \bM$ satisfies the hypotheses of the
  comparison lemma. Hence any sheaf on $\bM$ (in particular
  any manifold thought of as a representable sheaf) is determined
  by its values on $\Euc$.

  Hence one can think of a manifold as a certain kind of sheaf on
  the Euclidean site $\Euc$. This is an opposite perspective to
  the $C^{*}$-algebra approach, for example; there one describes a smooth
  manifold $M$ by its algebra of smooth functions, i.e. maps from
  $M$ to the fixed space $\R$. In the sheaf/stack picture
  one describes a manifold by the maps from the fixed spaces $R^{n}$
  to $M$.
  
  In fact if we restrict to manifolds of a fixed dimension $n$, 
  it is not hard to see that we need only $\mathrm{R}_{n}$, the
  site whose sole object is $\R^{n}$ and whose maps are
  the smooth self-maps $\R^{n} \ra \R^{n}$. However this is not
  a tremendous simplification and makes it hard to discuss manifolds
  of arbitrary (or nonconstant) dimension.
  
  What remains to be done is to characterize those sheaves
  which arise from smooth manifolds. They should be sheaves that
  are locally Euclidean in some sense. We will explain this in detail
  in Section~\ref{subsec:manifoldassheaf}.
\end{example}

From this example, which deals with sheaves (i.e.
stacks whose objects have no automorphisms), we turn to the
canonical example of a stack which is not a sheaf.

\begin{example}\label{ex:prinasclassifyingstack}
  Let $Z$ be a space and let 
  $G$ be a topological group and let
  $\Prin_{G}$ be the stack of principal $G$-bundles. 
  Then by Lemma \ref{lem:mapfromrepresentable}, 
  a map from $Z$ to $\Prin_{G}$ is just an element of
  the fiber $(\Prin_{G})_{Z}$, i.e. a principal bundle
  on $Z$. Hence we call $\Prin_{G}$ the \textit{classifying stack} of 
  $G$, in analogy to the classifying space in homotopy theory.
  We will often write it as $\BB G = \Prin_{G}$.
  
  Note the pros and cons of $\BB G$ versus the usual classifying space
  $BG$. The space $BG$ takes some work to construct, and it is only
  a classifying space once we identify homotopic maps; but at least it
  is an honest space. On the other hand, the stack 
  $\BB G$ is easy, even tautological, to construct, and it classifies
  bundles using (the natural stack extension of) ordinary maps.
  However this is at the expense of extending the category of
  spaces to the 2-category of stacks. So the tradeoff is between
  taking a quotient category (by homotopy) or embedding into a
  larger (2-)category (stacks).
\end{example}

We next define monomorphisms and epimorphisms of stacks.
Monomorphisms are uncomplicated.
\begin{Def}\label{def:monoepi}
  Let $a: F' \ra F$ be a map of fibered categories over $\CC$. 
  We say $a$ is a \textbf{monomorphism}
  if, for every object $C$ of $\CC$, the functor $a_{C}: F'_{C} \ra F_{C}$
  on fibers is fully faithful. 
\end{Def}

As in the case of sheaves and presheaves, there are different notions
of epimorphism for fibered categories versus stacks. We are most
interested in the following.
\begin{Def}
  We say a map of fibered categories 
  $a:F' \ra F$ over a site $(\CC,J)$
  is \textbf{covering} (French ``couvrant'') if, for every object
  $C$ of $\CC$, and every object $D$ of $F_{C}$,
  there is a covering family $\{f_{i}:C_{i} \ra C \}$
  and for every $i$ an object $D_{i}$ of $F'_{C_{i}}$ such that
  $a(D_{i}) \iso D|_{C_{i}}$.
  
  If $F'$ and $F$ are stacks, then we refer to a covering map 
  as an \textbf{epimorphism}.
\end{Def}

Another useful notion is that of the pullback of two maps of
stacks. It is not simply the 1-categorical notion of pullback;
this is one of the places where the 2-categorical nature of stacks 
appears explicitly.

\begin{Def}\label{def:pullbackoffibered}
  Let $F: \DD \ra \CC$, $F':\DD' \ra \CC$ and $F'':\DD'' \ra \CC$
  be three fibered categories and let $a: F' \ra F$ and $b:F'' \ra F$
  be maps of fibered categories. Then the \textbf{pullback} $F' \times_{F} F''$
  of the diagram 
  \[
    \xymatrix@1{ F' \ar[r]^{a} & F & F'' \ar[l]_{b}  }
  \]
  is defined as follows. The fiber $F_{C}$ over an object $C$ is
  the set of triples
  \[
    \{(D',D'',g) \: | \: D' \in \Ob F'_{C}, \: D'' \in \Ob F''_{C}, \:
                         (g: a(D') \ra b(D'')) \in \Ar F_{C} \}.
  \]
  An arrow from $(D_{1}',D_{1}'',g_{1})$ to $(D_{2}',D_{2}'',g_{2})$ 
  is a pair 
  \[
    (f':D_{1}' \ra D_{2}', f'':D_{1}'' \ra D_{2}'')
  \]
  with $f' \in \Ar F'_{C}$, $f'' \in \Ar F''_{C}$,
  such that
  \[
    g_{2} a(f') = b(f'') g_{1}.
  \]
  Finally, pullbacks are defined by, for $h:C' \ra C$ in $\CC$, 
  \begin{align*}
    h^{*}(D',D'',g) &= (h^{*}D',h^{*}D'',h^{*}g); \\
    h^{*}(f',f'')   &= (h^{*}f',h^{*}f''). \\
  \end{align*}
\end{Def}

We leave it to the reader to verify that the pullback of a
diagram of stacks is again a stack.

We now discuss the process of stackification, i.e. building
a stack out of a prestack.

\begin{Def}\label{def:stackification}
  Let $F: \DD \ra \CC$ be a fibered category. Then a 
  \textbf{stack associated to $F$}, or a \textbf{stackification}
  of $F$, is a stack $\hat{F}$ and a map of stacks  
  $i:F \ra \hat{F}$, such that for every stack $G$,
  composition with the map $i$ induces an equivalence of categories
  \[
    \HOM(\hat{F},G) \stackrel{i^{*}}{\ra} \HOM(F,G).
  \]
\end{Def}

\begin{proposition}\label{prop:existenceofstackification}
  Given a prestack $F: \DD \ra \CC$, there is a stack 
  associated to $F$, $i:F \ra \hat{F}$. It is unique up
  to a unique equivalence of stacks. The map $i$ is a covering
  monomorphism.
\end{proposition}
\textit{Remark.} In fact we do not really need to assume that $F$ is
a prestack; an associated stack exists for 
any fibered category \cite{Giraud:Book}.
However we will not need this result.

We will forgo the proof of Proposition~\ref{prop:existenceofstackification}, 
but we will describe the construction of $\hat{F}$. 
See \cite{LMB} or \cite{Giraud:Book} for more details. We will
describe it first in the language of sieves, as that is more 
elegant in this case. 
Given an object $C$ of a site $(\CC,J)$, we associate to it the category
whose objects are pairs 
\[
  \{(S,A) \: | \: S \in J(C) \text{~and~} A \in \Desc(S,F) \}.
\]
A map from $(S,A)$ to $(T,B)$ (over the same object $C$) is a 
map of descent data from $A|_{S \intersect T}$ to $B|_{S \intersect T}$.
Two such maps are identified if they agree on a further refinement.
The composition is evident (one has to restrict to a triple intersection
of sieves).

If one wants to use covering families instead of sieves, as
we typically do in this paper, the associated stack is defined
in the following way.
An object of the associated stack $\hat{F}$ over $C$
is given by the following data:
\begin{itemize}
\item a cover $\{C_{\alpha}\}$ of $C$; 
\item for each $C_{\alpha}$, an element $a_{\alpha} \in F_{C_{\alpha}}$;
\item for each $C_{\alpha \beta}$,
      an arrow $\psi_{\alpha \beta}: a_{\beta} \ra a_{\alpha}$, with 
      $\psi_{\alpha \beta} \cdot \psi_{\beta \gamma} = \psi_{\alpha \gamma}$
        on triple pullbacks.
\end{itemize}  
An arrow over $C$ from $(\{C_{\alpha}\},a,\psi)$
to $(C'_{\beta'},a',\psi)$ is a collection of
arrows 
\[
  \phi_{\alpha \beta'}: a_{\alpha}|_{C_{\alpha \beta'}}
    \ra a'_{\beta'}|_{C_{\alpha \beta'}}
\]
(where $C_{\alpha \beta'} = C_{\alpha} \times_{C} C_{\beta'}$ as usual),
compatible with the maps $\psi,\psi'$ in the obvious way.

Since stackification is a universal construction, we also
have a stackification of maps of prestacks in the usual way.
The following lemma, which generalizes a well-known fact in 
the case of sheaves, will be useful in Section~\ref{subsec:locallyrep}.
\begin{lemma}\label{lem:stackifyisiso}
  Let $a: F \ra G$ be a map of prestacks over the site $(\CC,J)$. 
  Then the stackified map
  $\hat{a}: \hat{F} \ra \hat{G}$ is an equivalence of stacks
  if and only if $a$ is a monomorphism and is covering.
\end{lemma}
\begin{proof}[Sketch of Proof.]
  Let $a:F \ra G$ be a covering monomorphism. We want a map
  $b: \hat{G} \ra \hat{F}$ which will be an inverse equivalence 
  to $\hat{a}$. By the definition of the associated stack, it is enough to
  construct a map $c:G \ra \hat{F}$. So, let $C \in \CC$ and let 
  $x \in G_{C}$. Since $a$ is covering, there is a cover $\{C_{\alpha} \ra C\}$
  and $y_{\alpha} \in F(C_{\alpha})$ with $a(y_{\alpha}) \iso x|_{C_{\alpha}}$.
  Since $a$ is a monomorphism (so $a:F_{C} \ra G_{C}$ is fully faithful)
  it is easy to show that $y_{\alpha}|_{C_{\alpha \beta}} \iso y_{\beta}|_{C_{\alpha \beta}}$
  and that these isomorphisms are coherent, i.e. they form descent data
  for $F$ over the cover $\{C_{\alpha} \ra C\}$. Let $z_{\alpha} = i(y_{\alpha})$.
  The descent data given by the $y_{\alpha}$ becomes descent data for the
  $z_{\alpha}$, and since $\hat{F}$ is a stack, we obtain an element
  $z \in \hat{F}(C)$. We define $c(x) = z$. This defines $c:G \ra \hat{F}$;
  it is straightforward to show that it is well-defined (and to define
  it on arrows as well as objects). The universal property of $\hat{G}$
  then gives a corresponding map $b: \hat{G} \ra \hat{F}$ such that
  $bj \iso c$.

  Note that for each $\alpha$,
  \[
    \hat{a}(c(x))|_{C_{\alpha}} = \hat{a} i(y_{\alpha}) 
      \iso j(a(y_{\alpha})) \iso j(x|_{C_{\alpha}};
  \]
  since $\hat{G}$ is a stack, this means $\hat{a}c \iso j$,
  so $\hat{a}bj \iso j$. Since $j$ is a monomorphism, this easily
  implies that $\hat{a}b \iso 1$. Also, $b\hat{a}i \iso bja \iso ca$,
  and $ca(y)|_{C_{\alpha}} \iso i(y|_{C_{\alpha}}$, so
  $ca \iso i$, which implies $b\hat{a}i \iso i$. Since $i$ is a monomorphism,
  we have $b\hat{a} \iso 1$.
\end{proof}

\section{Presentations of Stacks by Groupoids}\label{sec:pres}

In this section we will generalize Examples \ref{ex:sheafasstack} 
and \ref{ex:prinstackagain}  and see that we get
a large, flexible class of stacks. In particular we will
see how an orbifold defines a smooth stack.

We will work primarily in the topological category $\bT$ and the
differential category $\bM$ for definiteness. However
many of the results are quite general.

First we recall the definition of a \textit{topological groupoid},
that is, a groupoid object in the category of topological spaces.
\begin{Def}\label{def:topgroupoid}
  A \textbf{topological groupoid} $\GG$ is a pair $(\GGo, \GGl)$
  of topological spaces, together with five 
  continuous maps $s,t: \GG_{1} \ra \GG_{0}$, $u: \GG_{0} \ra \GG_{1}$, 
  $i: \GG_{1} \ra \GG_{1}$, 
  $m: \GG_{1} \times_{s,\GG_{0},t} \GG_{1} \ra \GG_{1}$,
  satisfying the relations given below.
  We call $\GG_{0}$ the space of \textbf{objects} and $\GG_{1}$
  the space of \textbf{arrows}.
  We use the notation $g: x \ra y$ to denote an arrow $g \in \GG_{1}$ with 
  $s(g) = x, t(g) = y$. Denoting $m(g,h)$ by $gh$, $i(g)$ by $g^{-1}$, 
  and $u(x)$ by $1_{x}$,
  the axioms are the familiar ones for a groupoid: for appropriate 
  $x, y \in \GG_{0}$, $g,h,k \in \GG_{1}$, 
  \begin{gather}\label{groupoidstructure}
    s(1_{x}) = t(1_{x}) = x, \quad s(g^{-1}) = t(g), \quad t(g^{-1}) = s(g); \\ 
    s(gh) = s(h), \quad t(gh) = t(g), \quad g(hk) = (gh)k; \\
    g 1_{x} = 1_{y} g = g, \quad  
    g^{-1} g = 1_{x}, \quad g g^{-1} = 1_{y} \quad \text{(for $g: x \ra y$).}
  \end{gather}
  
  If $s$ (hence $t$) is a local homeomorphism 
  we call $\GG$ an \textit{\'etale groupoid}. 
  If the map $(s,t): \GG_{1} \ra \GG_{0} \times \GG_{0}$ is
  proper, we call the groupoid \textit{proper}.
\end{Def}
Note: It is easy to show that $i$ is an involution of $\GG_{1}$; hence 
any property holds for $s$ if and only if it holds for $t$.

\begin{example}\label{ex:groupoids}
We have the following three simple examples of topological groupoids.
\begin{enumerate}
\item Given a space $X$, we can define the trivial groupoid, 
      also denoted $X$,
      with $X_{0} = X$ and with only identity arrows in $X_{1}$.
\item Given a topological group $K$ we can define a corresponding groupoid 
      $\KK$ with only one object, and with $K$ as the space of arrows.
\item Given a space $X$ and a right action of the group $K$ on $X$, 
      we can form the \textit{translation groupoid} 
      (or \textit{semi-direct product}) $X \rtimes K$ by
      $(X \rtimes K)_{0} = X$, $(X \rtimes K)_{1} = X \times K$, 
      $s(x,k) = xk$, $t(x,k) = x$,
      $u(x) = (x,1)$, $i(x,k) = (xk,k^{-1})$, 
      and $(x,k) \circ (y,h) = (x,kh)$.
\end{enumerate}
\end{example}

Note that the first example is $M \rtimes (e)$, 
and the second example is $\ast \rtimes K$. 

A \textit{strict homomorphism} $\phi: \GG \ra \HHH$ 
of topological groupoids is just a continuous functor, i.e. it is 
given by a pair of continuous maps $\phi_{0}: \GG_{0} \ra \HHH_{0}$
and $\phi_{1}: \GG_{1} \ra \HHH_{1}$ commuting with the structure maps of
$\GG$ and $\HHH$. We get a corresponding notion of strict
isomorphism. However one can also use a different notion of morphism,
explained in Section~\ref{subsec:moritaHS} below.

Given any topological groupoid $\GG$, one can
form its topological, or coarse, quotient 
$\GG_{\mathrm{top}} = \GG_{0}/\GG$ 
by identifying any two objects which have an arrow between them.
However this is not a very fine invariant of $\GG$ and is often not
a nice topological space.
For example, any groupoid that is categorically connected (i.e.
any two objects can be joined by an arrow) has
a topological quotient that is just a single point.
Note that if we take the topological 
quotient of a translation groupoid $X \rtimes K$ by identifying
objects with arrows between them, we get the topological quotient
$X/K$. 

We will see that the correct ``quotient'' of $\GG$ is actually
a \textit{stack} associated to $\GG$, which we will construct
below.

The definitions in the differential case are quite similar.
A \textit{smooth groupoid} is a topological groupoid $\GG$
where $\GGo,\GGl$ are smooth manifolds and all of the structure
maps are smooth, and additionally $s$ (hence $t$) is assumed to
be a submersion. (In particular this guarantees that $\GGl \times_{\GGo} \GGl$
is a smooth manifold.) For simplicity we will assume that both
$\GGo$ and $\GGl$ are in our category $\bM$, in particular that
they are Hausdorff. For some applications this is too restrictive,
but for the case of orbifolds, in particular, it is sufficient.

A smooth groupoid is \textit{\'etale} if $s$ is a local diffeomorphism.
It is \textit{proper} if it is proper as a topological groupoid.
A strict homomorphism of groupoids is defined as above but with
the maps required to be smooth.

\begin{example}\label{ex:smoothgroupoids}
We have examples of smooth groupoids paralleling those in 
Example~\ref{ex:groupoids}, and one new one.
\begin{enumerate}
\item Given a smooth manifold $M$ with a smooth right action of a 
      Lie group $K$, one can form the translation groupoid $M \rtimes K$
      as in the topological case. This includes the case 
      $M \rtimes (e)$, denoted simply by $M$, and $\ast \rtimes K = \KK$
      as above.
\item In particular, if $\KK$ is compact and 
      the action of $\KK$ is locally free,
      we obtain a smooth 
      groupoid that is easily seen to be \'etale and proper.
      Note that the ordinary quotient of $M$ by $\KK$ in this case
      has the structure of an orbifold.
\item Motivated by the previous example, one can define an 
      \textit{orbifold groupoid} to be a smooth, proper, \'etale groupoid
      (see \cite{Moerdijk:OrbGrpIntro} for an equivalent definition). 
      In fact this serves as a good alternate
      definition of an orbifold; in particular the correct notion
      of a map of orbifolds appears naturally in this language.
      We will see how this is related to thinking of orbifolds as
      stacks in Section~\ref{subsec:orbifoldasstack}.
\end{enumerate}
\end{example}

\subsection{The Stack Associated to a Groupoid}\label{subsec:assocstack}

Let $\GG$ be a topological groupoid in $\bT$. There is a natural way to
associate a fibered category over $\bT$ to $\GG$: define it 
as the presheaf of groupoids
\begin{align}\label{assocprestackeq}
  \Ob F_{X} &= \Maps(X,\GGo) \\
  \Ar F_{X} &= \Maps(X,\GGl) \label{assocprestackeqar}
\end{align}
with pullback functor induced by composition: for $f:X' \ra X$,
\begin{align}\label{assocprestackpullback}
  f^{*}(X \stackrel{a}{\ra} \GGo) &= X' \stackrel{af}{\ra} \GGo \\
  f^{*}(X \stackrel{g}{\ra} \GGl) &= X' \stackrel{gf}{\ra} \GGl. 
\end{align}
(Note that this is a \textit{strict} presheaf of groupoids---pullbacks
commute on the nose---which is somewhat convenient.)

It is easy to show that this fibered category is a prestack.
(This follows from the gluing axiom for maps of spaces, applied
to maps into $\GGl$.) However it is generally not a stack.
We have to \textit{stackify} it as in Section~\ref{sec:stackmaps}.
We define the stack associated to a topological groupoid
as $\uGG = \hat{F}$, where $F$ is the prestack 
associated to $\GG$ as in (\ref{assocprestackeq}),(\ref{assocprestackeqar}).
Hence a map from a space $A$ into $\uGG$ is given
locally by maps from open sets $U_{\alpha} \subset A$ into $\GGo$, with gluings given
by maps from $U_{\alpha \beta}$ into $\GGl$, satisfying the cocycle condition.

Note that we get a canonical map $p: \GGo \ra \uGG$; for,
$F(\GGo) = \Maps(\GGo,\GGo)$ contains the identity map,
giving a canonical map $\GGo \ra F$.
Composing with the canonical map $F \ra \hat{F}$ gives $p$.
Note that this map is an epimorphism. For, the canonical map
$F \ra \hat{F}$ is an epimorphism, and the map
$\GGo \ra F$ is just the identity map on objects. 
 
It is easy to see that this proces defines a functor 
\[
  S: (\text{Topological Groupoids}) \ra \St(\bT), \quad \GG \mapsto \uGG
\]
One can characterize this functor in another way, by looking
at an alternate definition of a map of groupoids, to which we now turn.

\subsection{Morita Equivalence and Hilsum-Skandalis Maps}
\label{subsec:moritaHS}

In this paper, we are not particularly interested in groupoids in themselves,
but rather as fine models for quotient ``spaces''. So the notion
of strict homomorphism is really too strong. In particular
one can have many strictly nonisomorphic groupoids with the
same associated stack.

\begin{example}\label{ex:chartcover}
Let $M$ be a smooth manifold and let 
$\UU = \{U_{i}\}$ be a cover of $M$ by charts. Let
$\GGo = \bigcup_{i} U_{i}$ and let $\GGl = \GGo \times_{M} \GGo$.
(So $\GGl$ is the disjoint union of the overlaps of the charts.) 
There is an obvious groupoid structure on $(\GGo,\GGl)$,
and an obvious strict homomorphism $p:\GG \ra M$.
Hence there is a map of the associated stacks 
$\underline{p}: \uGG \ra \underline{M}$. This map is an
equivalence of stacks. For, consider the prestack $F = \Maps(-,\GG)$
and the map $p: F \ra \underline{M}$. 
This map is easily
seen to be a monomorphism. It is also an epimorphism, since
every smooth map $f: N \ra M$ can locally be lifted to the charts $\GGo$.
Hence the induced map $\underline{p}$ on stacks is an equivalence.
So we would like to think of $M$ and $\GG$ as equivalent groupoids.
However the map $p$ of groupoids is, in general, 
clearly not a strict isomorphism, as it does not usually even have a section.
\end{example}

The map $p$ is an example of an \textit{essential equivalence}, 
defined as follows.
\begin{Def}\label{def:essequiv}
  A strict homomorphism $\phi: \GG \ra \HHH$ of topological groupoids 
  is an \textbf{essential equivalence} if both of the following
  conditions are satisfied:
  \begin{enumerate}
  \item \label{esssurj} $\phi$ is (topologically) 
        \textbf{essentially surjective:}
              the map 
              \[
                 s \circ \mathrm{pr}_{2}:\GGo \times_{\phi,\HHHo,t} 
                                                       \HHHl \ra \HHHo
              \]
              is a surjection admitting local sections;
  \item \label{fullyfaithful} $\phi$ is (topologically) 
        \textbf{fully faithful:} 
              the square
              \[
                \xymatrix{\GGl \ar[r]^-{\phi} \ar[d]_-{(s,t)}
                          & \HHHl \ar[d]^-{(s,t)} \\ 
                          \GGo \times \GGo \ar[r]^-{\phi \times \phi}    
                          & \HHHo \times \HHHo
                }
              \]
              is a fiber product.
  \end{enumerate}
\end{Def}
\begin{remark}\label{rem:etaleissimplersurjection}
  If we restrict to the category of \'etale groupoids,
  then to ensure $\phi$ is topologically essentially surjective,
  it is enough to require that $\phi$ is an open surjection.

  If we work in the category of smooth groupoids, we require in
  (\ref{esssurj}) that $\phi$ be a surjective submersion.
\end{remark}

It is important to remember that an essential equivalence does
not generally have an inverse. If we forget the topology,
an essential equivalence is a categorical equivalence, so it
has a quasi-inverse, but this will usually be discontinuous.
However essential equivalences become honest equivalences when
we pass to stacks. See \cite{Pronk:Bicategories} for a detailed exposition of 
the following.
\begin{proposition}\label{prop:essequivgivesstackiso}
  Let $\phi: \GG \ra \HHH$ be an strict homomorphism of 
  topological groupoids. Then the associated map of stacks
  $\underline{\phi}:\GG \ra \HHH$ is an equivalence if and only if
  $\phi$ is an essential equivalence. 
\end{proposition}
\begin{proof}
First, suppose that $\phi$ is an essential equivalence.
Let $F_{\GG} = \Maps(-,\GG)$ and $F_{\HHH} = \Maps(-,\HHH)$ be the
prestacks represented by $\GG,\HHH$. We need to show that the map 
$F_{\phi}:F_{\GG} \ra F_{\HHH}$ 
is a monomorphism and an epimorphism. It is easy to see that
since $\phi$ is topologically fully faithful, the functor
$F_{\phi}: F_{\GG} \ra F_{\HHH}$ is fully faithful on each fiber
$(F_{\GG})_{X}$, for any space $X$.

Now let $X$ be a space and let $f: X \ra \HHHo$ be an object
of $(F_{\HHH})_{X}$. Since $\phi$ is topologically essentially surjective,
there is an open cover $\{U_{i}\}$ of $\HHHo$ and sections 
$\sigma_{i}:U_{i} \ra \GGo \times_{\HHHo} \HHHl$ of the map $s \pr_{2}$.
Let $V_{i} = f^{-1}(U_{i})$, an open cover of $X$, and consider
the maps $g_{i}: V_{i} \ra \GGo$ and $h_{i}:V_{i} \ra \HHHl$ given by
\[
  g_{i} = \pr_{1} \sigma_{i} f, \quad h_{i} = \pr_{2} \sigma_{i} f.
\]
Then $g_{i}$ is an object of $(F_{\GG})_{V_{i}}$,
and $h_{i}$ is an isomorphism in $(F_{\HHH})_{V_{i}}$ between
$F_{\phi}(g_{i})$ and $f|_{V_{i}}$. Hence $F_{\phi}$ is an epimorphism,
and the associated map of stacks $\underline{\phi}$ is an equivalence.

Now assume that $\underline{\phi}$ is an equivalence of stacks.
By Lemma~\ref{lem:stackifyisiso}, the map $F_{\phi}$ of prestacks
must be a monomorphism and an epimorphism. Once again it is easy
to verify that $\phi$ must be topologically fully faithful.
To show that $\phi$ must be topologically essentially surjective,
apply the definition of an epimorphism of stacks to the
space $\HHHo$. This gives an open cover of $\HHHo$ and
sections of $s \pr_{2}$ over the sets of this cover.
\end{proof}

Motivated by this result, instead of passing to stacks, 
one can work in the category of groupoids, and 
formally invert the essential equivalences. (I.e. one 
\textit{localizes} the category at the essential equivalences.)
More explicitly, following \cite{MoerdijkPronk:OrbifoldsSheavesGroupoids}, we say 
a \textit{generalized map} from $\GG$ to $\HHH$
is an equivalence class of diagrams of the form
\[
  \xymatrix@1{ \GG & \GG' \ar[l]_{\varepsilon} \ar[r]^{\phi} & \HHH}
\]
where $\varepsilon$ is an equivalence. Another diagram
\[
  \xymatrix@1{ \GG & \GG'' \ar[l]_{\delta} \ar[r]^{\psi} & \HHH}
\]
is equivalent to the first if there is a homomorphism
$\gamma: G'' \ra G'$ with $\varepsilon \gamma \iso \delta$ and
$\phi \gamma \iso \psi$. 

One can show that the resulting category is equivalent to a full
subcategory of the category of topological stacks.
(In fact there is a \textit{bicategory} structure on
groupoids, and this leads on localization to the 2-category structure
on stacks. See Remark~\ref{rem:HSbicategory} below, and for a very
full treatment, see Pronk \cite{Pronk:Bicategories}.)

There is another way to describe the resulting localized category
which is pleasantly concrete, due to Hilsum and Skandalis. It
also relates nicely to the $C^{*}$-algebra treatment of groupoids.
We present it briefly.

Let $\GG$ be a topological 
groupoid. A \textit{right action} of $\GG$ on a space
$X$ consists of two continuous maps $\pi: X \ra \GG_{0}$ (often called 
the ``moment map''; we will call it the ``base map'' to avoid
confusion with symplectic geometry usage)
and 
\begin{align*}
  m: X \times_{\pi,\GG_{0},t} \GG_{1} = \{ (x,g) \: | \: \pi(x) = t(g) \} 
           &\ra X \\
     (x,g) &\mapsto xg    
\end{align*}
(the action map) such that 
\[
  (xg)h = x(gh), \quad x1 = x, \quad \pi(xg) = s(g).
\]

Given a right action of $\GG$ on $X$ we can form the semidirect
product, or translation, groupoid, generalizing the third example
from the previous section: we define
\[
  (X \rtimes \GG)_{0} = X, \qquad 
    (X \rtimes \GG)_{1} = X \times_{\pi,\GG_{0},t} \GG_{1} 
\]
and $s(x,g) = xg$, $t(x,g) = x$, $(x,g)(xg,h) = (x,gh)$,
$i(x,g) = (xg,g^{-1})$.
Note that $\pi$ extends as a covariant functor $\pi: X \rtimes \GG \ra \GG$
by defining $\pi(x,g) = g$.
Unless otherwise specified all actions will be from the right,
so by a \textit{$\GG$-space} we will mean a space with a 
given right $\GG$-action.

A left action of $\GG$ is a right action of the opposite groupoid
$\GG^{\mathrm{op}}$. 

\begin{Def}
  Given a space $B$, a \textbf{(right) $\GG$-bundle} over $B$ 
  is a (right) $\GG$-space
  $E$ and a continuous $\GG$-invariant map $p: E \ra B$ (so 
  $p(xg) = p(x)$). It is \textbf{principal} if $p$ is a 
  open surjection and
  \begin{equation}\label{prinbundle}
    E \times_{\GG_{0}} \GG_{1} \ra E \times_{B} E, \quad (e,g) \mapsto (e,eg)
  \end{equation}
  is a homeomorphism.

  In the smooth case, we require $p$ to be a surjective submersion
  and the map in (\ref{prinbundle}) to be a diffeomorphism.
\end{Def}
This reduces to the usual notion of a principal bundle in the
case of a topological (resp. Lie) group, i.e. when $\GG_{0}$ is a point.

\begin{Def}\label{def:HSmor}
  Let $\GG$ and $\HHH$ be topological groupoids. A Hilsum-Skandalis (HS)
  morphism $f = (P,\sigma,\tau):\GG \ra \HHH$ consists of a space $P$, 
  maps $\sigma : P \ra \GG_{0}$, $\tau : P \ra \HHH_{0}$, 
  a left action of $\GG$ on $P$
  with the base map $\sigma$, a right action of $\HHH$ on
  $P$ with base map $\tau $, such that:
  \begin{enumerate}
  \item $\sigma$ is $\HHH$-invariant, $\tau $ is $\GG$-invariant;
  \item the actions of $\GG$ and $\HHH$ on $P$ are compatible: given $p \in P$,
        $(gp)h = g(ph)$;
  \item $\sigma: P \ra \GG_{0}$, as an $\HHH$-bundle with moment map $\tau$,
        is principal. 
  \end{enumerate}
  
  Given two maps $f = (P,\sigma,\tau)$ and $f' = (P',\sigma',\tau')$ 
  from $\GG$ to $\HHH$, 
  a \textbf{2-isomorphism} from
  $f$ to $f'$ is a $\GG$- and $\HHH$-equivariant homeomorphism 
  $\phi: P \ra P'$ satisfying $\sigma' \phi = \sigma$ and $\tau' \phi = \tau$.

  One makes the obvious analogous definitions in the smooth case.
\end{Def}

The composition of two HS morphisms $P: \GG \ra \HHH$ and
$Q: \HHH \ra \KK$ is defined by dividing out $P \times_{\HHH_{0}} Q$
by the action of $\HHH$: $(p,q)h = (ph,h^{-1}q)$, and taking
the obvious actions of $\GG$ and $\KK$. The identity
map of $\GG$ is represented by the diagram
$\GGo \la \GGl \ra \GGo$ where the maps are the source and target
maps and the actions are the left and right multiplication.
A map $\GGo \la P \ra \HHHo$ is invertible (modulo the remark below) 
if both $s_{P}$ and $t_{P}$ are principal bundles.

\begin{remark}\label{rem:HSbicategory}
Since fiber products are only associative up to a (unique) natural
isomorphism, the same is true of composition of HS maps: they associate
only up to a (canonical, coherent) 2-isomorphism.
Similarly for identity maps and inverses. 
Hence the HS maps as defined do not form a category, but rather
a so-called \textit{bicategory}. This is not a very big worry.
It is analogous to the fact that a stack defines a lax presheaf
of groupoids and not an honest one. 
\end{remark}

Given a strict homomorphism $\phi: \GG \ra \HHH$, there is
an associated HS morphism given by
\[
  \xymatrix{ \GGo & \GGo \times_{\phi, \HHHo, t} \HHHl 
                      \ar[l]_-{\mathrm{pr}_{1}} 
                      \ar[r]^-{s \circ \mathrm{pr}_{2}} 
                  & \HHHo
  }
\]
where the actions of $\GG$ and $\HHH$ on $\GGo \times_{\HHHo} \HHHl$
are given by $g \cdot (x,h) \cdot h' = (t(g),\phi(g) h h')$.
It is easy to check that this defines an HS morphism.

It is also easy to check the following.
\begin{proposition}\label{prop:essequivgivesHSequiv}
  A strict homomorphism $\phi: \GG \ra \HHH$ is an
  essential equivalence if and only if the resulting HS map is an 
  HS equivalence.
\end{proposition}

\subsection{Locally Representable Stacks}\label{subsec:locallyrep}

In the algebraic setting, there is a standard way (\cite{DeligneMumford})
to characterize the stacks which one obtains
using the construction of the previous section, which shows
that these are a good generalization of schemes. We will explain how things
work in the topological and smooth cases; they are a little different,
because pullbacks do not always exist in the smooth case.

We begin with a few general preliminaries. We recall that a stack 
$F$ over $\CC$ is \textit{representable} if $F$ is equivalent to $\uC$
for some object $C$ of $\CC$. As usual, given an object 
$C$ of $\CC$, we will write $C$ instead of $\uC$ where
there will be no confusion.

\begin{Def}\label{def:repmap}
  Let $F, G$ be stacks over $\CC$. We say that a map $a:F \ra G$
  is \textbf{representable} if for every object $C$ of $\CC$,
  the pullback stack $C \times_{G} F$ is representable.
\end{Def}

\begin{Def}\label{def:locallyrep}
  Let $F$ be a stack over $\CC$. We say $F$ is \textbf{locally representable}
  if there is an object $C$ of $\CC$ and a representable epimorphism
  $p: C \ra F$.
\end{Def}
We will call such a $p:C \ra F$ a \textit{presentation} of the stack $F$.
It is also often called a \textit{chart} or \textit{atlas}. 
For, in the manifold
case, if we are given an atlas $\{U_{\alpha} \}$ for a manifold $M$,
then the canonical map $\coprod_{\alpha} U_{\alpha} \ra M$ is a 
presentation of $M$ by an object in $\Euc$.

We note that a map of topological spaces $f:X \ra Y$, considered as a
map of stacks, is an epimorphism if and only if it has local sections.
For, suppose that $f: X \ra Y$ is an epimorphism. Consider the
identity map $1_{Y}$ as an element of $\underline{Y}(Y)$.
Since $f$ is an epimorphism, there must be a cover $\{U_{\alpha}\}$of $Y$ and
elements $\sigma_{\alpha} \in \underline{X}(U_{\alpha})$ with
$f(\sigma_{\alpha}) = (1_{Y})|_{U\alpha}$. But this just says that
$\sigma_{\alpha}:U_{\alpha} \ra X$ is a local section.

Conversely, suppose that $f$ has local sections $\{\sigma_{\alpha}\}$. Given a map
$g: A \ra Y$, we can compose with the $\sigma_{\alpha}$ to get local lifts
of $g$ to $X$, showing that $f$ is an epimorphism of stacks.

This is an example of a general phenomenon: requiring something to be
true on the stack level means it has to be \textit{locally effective}.
For a map $f$ to be an epimorphism, it is not enough for $f$ 
to be surjective (i.e. to
have sections over single points); one has to have sections
over each set of some open cover.

Any property of a map in $\CC$ which is stable under pullback
can be made into a property of a representable map of stacks.
\begin{Def}\label{def:propertiesofrepmaps}
  Let $P$ be a property of maps in $\CC$ that is stable under
  pullback. We say that a representable map $f:F \ra G$
  \textbf{has property $P$} if, for every object $C$ of $\CC$ and
  map $g: C \ra G$, the projection $g^{*}f: C \times_{G} F \ra C$ has property $P$.
\end{Def}

This works quite well in the topological setting, to which
we will soon turn. (We will see that we need to be more careful in the smooth
setting below.) Def. \ref{def:propertiesofrepmaps} allows us to say when
a (representable) map of stacks over $\bT$ (or more generally over $\bTX$)
is open, \'etale, surjective, or an embedding.

We will want the following general lemma:
\begin{lemma}\label{lem:diagonalgivessheaf}
  Let $F$ be a stack over $\CC$ and let $\Delta:F \ra F \times F$ be the diagonal map. 
  Let $C$ be an object of $\CC$ and let $c: C \ra F \times F$ be a map. Then the
  pullback of stacks $G = F \times_{\Delta,F \times F,c} C$ is equivalent to a sheaf
  (that is, each object has trivial automorphisms).
  An object of $G_{A}$ is given by a map $a: A \ra C$
  and an arrow $\alpha: c^{1}a \ra c^{2}a$ in $F_{A}$.
\end{lemma}
\begin{proof}
  Let $A$ be an object of $\CC$.
  By definition of the pullback of stacks, an object of $G_{A}$ 
  is given by a map $a: A \ra C$, an object $f \in F_{C}$, and an arrow from
  $ca = (c^{1}a,c^{2}a)$ to $(f,f)$ in $(F \times F)_{A}$, in other words,
  a pair of arrows $\alpha^{1}:c^{1}a \ra f, \: \alpha^{2}:c^{2}a \ra f$.
  An arrow from $(a,f,\alpha^{1},\alpha^{2})$ to itself is given by 
  an arrow $\phi:f \ra f$ such that $\phi \alpha^{i} = \alpha^{i}$, $i=1,2$.
  (Usually one would have an arrow from $a$ to $a$ as well, but $a$ has
  no nontrivial automorphisms since $\uC$ is a sheaf.)
  But this requires that $\phi = 1$, so $(a,f,\alpha^{1},\alpha^{2})$ has
  no nontrivial automorphisms. This implies that each category
  $G_{A}$ is equivalent to a discrete category, so $G$ is equivalent to a sheaf.
  
  In fact we can make this equivalence explicit as follows. Given
  $(a,f,\alpha^{1},\alpha^{2})$, this is canonically isomorphic to 
  $(a,c^{1}a,1,(\alpha^{2})^{-1} \alpha^{1})$.
  So the latter is a canonical representative for the isomorphism class
  of objects represented by $(a,f,\alpha^{1},\alpha^{2})$.
  Letting $\alpha = (\alpha^{2})^{-1} \alpha^{1}$, we see that
  an object of $G_{A}$ (after passing to the equivalent sheaf)  
  can be represented as $(a,\alpha)$ as desired.
\end{proof}

The following is a simple diagram chase.
\begin{cor}\label{cor:pullbackbydiag} 
  Let $F$ be a stack over $\CC$ and let $\Delta:F \ra F \times F$ 
  be the diagonal map. 
  Let $C$ be an object of $\CC$ and let $p: C \ra F$ be a map. 
  Assume that either $\Delta$ or $p$ is representable. Then the
  pullback of stacks $F \times_{\Delta,F \times F,p \times p} (C \times C)$ 
  is equivalent to the pullback $C \times_{p,F,p} C$, and both
  are equivalent to sheaves.
\end{cor}

Given a groupoid $\GG$, a space $X$, and two maps $f,g:X \ra G_{1}$,
we write $f \cdot g$ to denote the product in the groupoid, that is,
\[
  f \cdot g = m \circ (f,g).
\]

\begin{proposition}\label{prop:presentationgivesgroupoid}
  Let $\GG$ be a topological groupoid (in $\TT$) and let $\uGG$ be the 
  associated stack as in Sec. \ref{subsec:assocstack}. 
  Then $\uGG$ is locally representable; explicitly, 
  the canonical epimorphism $p:\GGo \ra \uGG$ is representable, and in particular,  
  $\GGo \times_{\uGG} \GGo \iso \GGl$. In addition,
  the diagonal map $\Delta: \uGG \ra \uGG \times \uGG$ is also representable.

  Conversely, suppose that $F$ is a locally representable stack over $\bT$
  and $p:C_{0} \ra F$ is a presentation. 
  Let $F' = C_{0} \times_{F} C_{0}$. Then
  $F'$ is represented by some object $C_{1}$ of $\CC$; 
  the pair $\mathcal{C} = (C_{1},C_{0})$ 
  has a canonical structure of a groupoid
  object in $\CC$, where the source and target maps are
  the projections $\pi_{1},\pi_{2}$; and the stack $\underline{\mathcal{C}}$
  associated to $\mathcal{C}$ is canonically equivalent to $F$.
\end{proposition}
\begin{proof}
  First we show that $\Delta$ is representable. Let $X$ be a topological space in $\TT$
  and let $a:X \ra \uGG \times \uGG$ be a map. A map from $X$ to $\uGG \times \uGG$
  is the same thing as an element of $(\uGG \times \uGG)(X)$, which in turn is a pair
  $(a^{1},a^{2})$ of elements of $\uGG(X)$. These elements $a^{1},a^{2}$, by definition of
  the associated stack (see the discussion following 
  Prop. \ref{prop:existenceofstackification}), 
  are given by the following data:
  \begin{itemize}
  \item a cover $\{U^{i}_{\alpha}\}$ of $X$ ($i=1,2$);
  \item for each $U^{i}_{\alpha}$, a map $a^{i}_{\alpha}: U^{i}_{\alpha} \ra \GGo$;
  \item for each intersection $U^{i}_{\alpha \beta} = U^{i}_{\alpha} \intersect U^{i}_{\beta}$,
        a map $\psi^{i}_{\alpha \beta}: U^{i}_{\alpha \beta} \ra \GGl$, with
        $s \circ \psi^{i}_{\alpha \beta} = a^{i}_{\beta}|_{U^{i}_{\alpha \beta}}$, 
        $t \circ \psi^{i}_{\alpha \beta} = a^{i}_{\alpha}|_{U^{i}_{\alpha \beta}}$, and 
        $\psi^{i}_{\alpha \beta} \cdot \psi^{i}_{\beta \gamma} = \psi^{i}_{\alpha \gamma}$
        on triple intersections.
  \end{itemize}  
  By passing to a common refinement, we can assume that the covers are the same, 
  $U^{1}_{\alpha} = U^{2}_{\alpha} = U_{\alpha}$ for all $\alpha$.

  We first determine what a map from a space $Y$ to
  $X \times_{a,F \times F,\Delta} F$ should be, if the pullback exists as a space.
  By Lemma \ref{lem:diagonalgivessheaf}, a map from $Y$ is given by 
  a map $h:Y \ra X$ and a transformation $\phi$ from $a^{1}h$ to $a^{2}h$.
  By definition, this is in turn given by a map $\phi_{\alpha}:h^{-1}(U_{\alpha}) \ra G_{1}$
  for every $\alpha$, with $s \circ \phi_{\alpha} = a^{1}_{\alpha}h$,
  $t \circ \phi_{\alpha} = a^{2}_{\alpha}h$, and 
  \begin{equation}\label{psiphiequation}
    \psi^{2}_{\alpha \beta} \cdot \phi_{\beta} = \phi_{\alpha} \cdot \psi^{1}_{\alpha \beta} 
  \end{equation}
  on every overlap $h^{-1}(U_{\alpha \beta})$.
  
  Now we seek to construct the pullback $X \times_{a,F \times F,\Delta} F$ as a topological
  space over $X$. First we will construct it locally, over each $U_{\alpha}$. Define
  $A_{\alpha}$ as the following pullback:
  \[
    \xymatrix{ A_{\alpha} \ar[r] \ar[d] & \GGl \ar[d]^{(s,t)} \\
               U_{\alpha} \ar[r]_-{(a^{1}_{\alpha},a^{2}_{\alpha})} & \GGo \times \GGo.
    }
  \]
  A section of $A_{\alpha}$ over $U_{\alpha}$ is given by 
  $\phi_{\alpha}:U_{\alpha} \ra \GGl$ with $s \phi_{\alpha} = a^{1}_{\alpha}$,
  $t \phi_{\alpha} = a^{2}_{\alpha}$. 
  More generally, given a map $h:Y \ra X$, a lift of $h|_{h^{-1}(U_{\alpha})}$ to $A_{\alpha}$
  is a map $\phi_{\alpha}:h^{-1}(U_{\alpha}) \ra \GGl$ 
  with $s \phi_{\alpha} = a^{1}_{\alpha}h$,
  $t \phi_{\alpha} = a^{2}_{\alpha}h$. 

  Now we glue the spaces $A_{\alpha}$ together. Over each double intersection 
  $U_{\alpha \beta}$, define a transition map 
  $\rho_{\alpha \beta}: A_{\beta}|_{U_{\alpha \beta}} \ra A_{\alpha}|_{U_{\alpha \beta}}$
  by 
  \[
    \rho_{\alpha \beta}(c,g) 
     = (c, \psi^{2}_{\alpha \beta}(c) \cdot g \cdot \psi^{1}_{\beta \alpha}(c))
  \]
  where $c \in U_{\alpha \beta}$, 
  $s(g) = a^{1}_{\beta}(c)$, and $t(g) = a^{2}_{\beta}(g)$, 
  so $(c,g) \in A_{\beta}|_{U_{\alpha \beta}}$.

  On a triple intersection $U_{\alpha \beta \gamma}$, we have 
  \begin{align*}
    \rho_{\alpha \beta} \circ \rho_{\beta \gamma} (c,g) 
      &= (c, \psi^{2}_{\alpha \beta}(c) \cdot \psi^{2}_{\beta \gamma}(c) \cdot g 
              \cdot \psi^{1}_{\gamma \beta}(c) \cdot \psi^{1}_{\beta \alpha}(c))\\
      &= (c, \psi^{2}_{\alpha \gamma}(c) \cdot g 
              \cdot \psi^{1}_{\gamma \alpha}(c))\\
      &= \rho_{\alpha \gamma}(c,g).
  \end{align*}
  Hence we can glue the spaces $A_{\alpha}$ together using the homeomorphisms
  $\rho_{\alpha \beta}$. Call the resulting space $A$; we claim it is the
  desired pullback. For, a lift $\phi$ of $h:Y \ra X$ to $A$ is given by a family of maps
  $\sigma_{\alpha}: h^{-1}(U_{\alpha}) \ra A_{\alpha}$ which are compatible under the gluings,
  i.e. such that 
  \[
    \rho_{\alpha \beta} \circ \sigma_{\beta} = \sigma_{\alpha}
  \]
  on $U_{\alpha \beta}$. But each local lift to $A_{\alpha}$ is given by 
  $\phi_{\alpha}: U_{\alpha} \ra \GGl$, and the above equation translates to
  \[
    \psi^{2}_{\alpha \beta} \cdot \phi_{\beta} \cdot \psi^{1}_{\beta \alpha} 
        = \phi_{\alpha}   
  \]
  which is the same as (\ref{psiphiequation}). Hence $A$ is the desired pullback,
  showing that the map $\Delta$ is representable.

  We now show that $p:\GGo \ra \uGG$ is representable; the proof is quite 
  similar to the proof of the previous statement. (In fact
  one can show that the representability of $\Delta$ implies the representability
  of $p$, see \cite{LMB}. However we will show it explicitly, for comparison to the
  smooth case.) Consider a space $X$ and
  a map $a: X \ra \uGG$. Explicitly, $a$ is given by
  $\{(U_{\alpha},a_{\alpha}:U_{\alpha} \ra \GGo, \psi_{\alpha \beta}:U_{\alpha \beta} \ra \GGl)\}$
  as we saw above. Define $B_{\alpha}$ to be the following pullback:
  \[
    \xymatrix{ B_{\alpha} \ar[r] \ar[d] & \GGl \ar[d]^{t} \\
               U_{\alpha} \ar[r]_{a_{\alpha}} & \GGo. 
    }
  \]
  On an overlap $U_{\alpha \beta}$, we have a transition function $\rho_{\alpha \beta}$ defined
  by 
  \[
    \rho_{\alpha \beta}(x,g) = (x,\psi_{\alpha \beta}(x) \cdot g).
  \]
  As before, the cocycle condition on the $\psi_{\alpha \beta}$ implies that
  we can use the $\rho_{\alpha \beta}$ to glue the $B_{\alpha}$ together;
  call the resulting space $B$. We claim that this is the desired pullback.

  Let $h:Y \ra X$ be a map. Then a lift of $h$ to $B$ is given by a family of maps 
  $\phi_{\alpha}: h^{-1}(U_{\alpha}) \ra \GGl$ with $t \phi_{\alpha} = a_{\alpha}h$
  and $\psi_{\alpha \beta} \cdot \phi_{\beta} = \phi_{\alpha}$.
  Note that $s \phi_{\alpha}=s \phi_{\beta}$, so that $f = s \phi: Y \ra \GGo$
  is well-defined, and $\phi$ gives a map of descent data from $f$ (a global
  object, considered as a descent datum) to $a$. This is exactly the data of
  a map from $Y$ to the pullback $X \times_{a,\uGG,p} \GGo$.

  Last, we show that $\GGo \times_{\uGG} \GGo \iso \GGl$. This is quite
  simple. Suppose we are given a space $X$, maps $h,k:X \ra G_{0}$,
  and a transformation $\phi: ph \ra pk$ of maps from $X$ to $\uGG$.
  By definition, such a transformation is a map $\phi: X \ra \GGl$ such that
  $s \circ \phi =h$ and $t \circ \phi = k$. But this just says that the
  square
  \[
    \xymatrix{ \GGl \ar[r]^{t} \ar[d]_{s} & \GGo \ar[d]^{p} \\
               \GGo \ar[r]_{p} & \uGG. 
  }
  \]
  is a pullback, as desired.

  Now we show the converse. Assume that $F$ is a locally representable stack,
  with presentation $p:C_{0} \ra F$.
  That $F' = C_{0} \times_{F} C_{0}$ is representable follows directly from the fact
  that $p:C_{0} \ra F$ is representable. Let $C_{1}$ be the space representing $F'$.
  The unit map $u:C_{0} \ra C_{1}$
  is given by the diagonal map $C_{0} \ra C_{0} \times_{F} C_{0}$. 
  We need to identify the product map $m$ and the inverse map $i$. 
  The inverse map is simply the canonical map switching factors
  in $C_{0} \times_{F} C_{0}$.
  The product map is given by 
  \[
    C_{1} \times_{C_{0}} C_{1} 
     \iso (C_{0} \times_{F} C_{0}) \times_{C_{0}}
            (C_{0} \times_{F} C_{0})
     \iso C_{0} \times_{F} C_{0} \times_{F} C_{0}
     \stackrel{\pi_{1,3}}{\lra} C_{0} \times_{F} C_{0}.
  \]
  We need to be careful to justify the second isomorphism above,
  since each pullback over $F$ is a pullback in the stack sense.
  A map from a space $A$ into $C_{0} \times_{F} C_{0}$   is given by 
  \[
    (a:A \ra C_{0}, b:A \ra C_{0},\phi:pa \ra pb).
  \]
  Hence a map from $A$ into $(C_{0} \times_{F} C_{0}) \times_{C_{0}} (C_{0} \times_{F} C_{0})$
  is given by 
  \[
   (a:A \ra C_{0}, b:A \ra C_{0},c: A \ra C_{0},\phi:pa \ra pb,\psi:pb \ra pc).
  \]
  A map from $A$ into $C_{0} \times_{F} C_{0} \times_{F} C_{0}$ is given by 
  \[
    (a:A \ra C_{0}, b:A \ra C_{0},c: A \ra C_{0},\phi:pa \ra pb,\psi:pb \ra pc,
      \chi:pa \ra pc)
  \]  
  where we require $\chi = \psi \phi$. These data clearly determine each other.
  Hence the pullbacks are equivalent as stacks, and since they are representable,
  they are isomorphic as spaces.

  It is straightforward to verify that these maps satisfy the axioms for a groupoid.
  For example, we will verify the identity $m \circ (us,1) = 1$, i.e. $u(g) \cdot g = g$
  for all $g \in C_{1}$. By our definitions of $m$, $u = \Delta$ and $s= \pi_{1}$, for
  any space $A$ and map $f: A \ra C_{1} = C_{0} \times_{F} C_{0}$ given by 
  $(f_{1},f_{2},\phi: pf_{1} \ra pf_{2})$, we have
  \[
    m \circ (us,1) \circ f = m \circ (f_{1},f_{1},f_{1},f_{2}) = (f_{1},f_{2}) = f.
  \]

  To construct the stack $\underline{\mathcal{C}}$,
  we first construct the prestack $G$ represented by $\mathcal{C}$
  as in Section \ref{subsec:assocstack}.
  There is a canonical map of stacks $j: G \ra F$ given as follows.
  Let $U$ be  an object of $\CC$ and let $(\arr{U}{a}{C_{0}}) \in \Ob G_{U}$.
  On objects, the map $j$ is defined by
  \[
    j_{U}(a) = pa 
  \]
  where as usual we identify the map $pa:\underline{U} \ra F$
  with an element of $F(U)$. As for arrows, consider
  $(\arr{U}{g}{C_{1}}) \in \Ar G_{U}$. By the definition
  of the pullback of stacks (\ref{def:pullbackoffibered}), $g$ is given by 
  a triple $(b,c,\gamma)$ with $b,c: U \ra C_{0}$ and
  $(\arr{pb}{\gamma}{pc}) \in \Ar F(U)$. 
  Then $j_{U}(g)$ is defined by
  \[
    j_{U}(g) = j_{U}(b,c,\gamma) = \gamma.
  \]
  Compatibility with the pullback maps is easy to check.

  The map $j$ is both a monomorphism and an epimorphism.
  For, each $j_{U}$ is clearly fully faithful, so $j$ is a monomorphism. 
  Also, since $j$ is just defined by $p$ on objects, and
  $p$ is an epimorphism, we see that $j$ is an epimorphism.

  Now let $\underline{\mathcal{C}}$ be the stack associated to
  $G$. By the universal property of the associated stack,
  we get an induced map $k:\underline{\mathcal{C}} \ra F$.
  Using Lemma~\ref{lem:stackifyisiso}, we see that $k$ is
  an equivalence.
\end{proof}

Hence we see that the stacks arising from groupoids as in 
Section~\ref{subsec:assocstack} are exactly the locally
representable stacks. We note that this statement can be sharpened
into an equivalence of bicategories; see \cite{Pronk:Bicategories}.

In the smooth case we must be a bit more careful. 
First, since we want to be able to take (transverse) pullbacks,
we will consider stacks over $\bM$ instead of over $\Euc$ (the
``pragmatic'' perspective).

Note that since pullbacks of manifolds do not generally exist, the condition
that a map be representable is quite strong. Note that if $M, N$ are
smooth manifolds and $f:M \ra N$ is a smooth map, then
$f$ is representable (when considered as a map between stacks on $\bM$)
if and only if every pullback of the form $M \times_{f,N,a} A$
exists, where $A \in \bM$ and $a:A \ra N$ is smooth.
But note the following lemma.
\begin{lemma}\label{lem:smoothrepissubmersion}
  Let $f:M \ra N$ be a smooth map. A pullback $M \times_{f,N,a} A$ exists as a smooth manifold
  for every smooth map $a: A \ra N$ if and only if the map $f$ is a
  submersion.
\end{lemma}
\begin{proof}
  If $f$ is a submersion, the pullback exists by standard trasnversality arguments,
  e.g. \cite{GuilleminPollack}. Conversely, suppose that $f$ is not a submersion. 
  First we show that we can reduce to the case where $N = \R$. There is some
  $p \in M$ and some $v \in T_{f(p)}N \setminus f_{*}(T_{p}M)$.
  Let $c:\R \ra N$ be a regular curve with $c(0)=p$, $c'(0)=v$. If the pullback
  $M \times_{f,N,c} \R$ does not exist, we are done. So assume it does exist; then
  the map $c^{*}f: M \times_{N} \R \ra \R$ is a smooth real-valued function with
  a critical point at $p$. 
 
  So from now one assume that $N = \R$ and $p$ is a critical point of the
  real-valued function $f:M \ra \R$. We can clearly assume that $f(p)=0$,
  and working locally on $M$, we can assume that $M = U \subset \R^{m}$
  and $p = 0$. Now let $A = \R$ and $a:\R \ra \R$ be $a(y)=y^{2}$.
  Then the pullback is (as a set) 
  \[
    M \times_{N} A = \{(x,y) \in R^{n} \times \R \: | \: y^{2} = f(x) \}.
  \]
  This has a singularity at the origin, since $0$ is a critical point of $f$.
  Hence $M \times_{N} A$ does not exist as a smooth manifold.
\end{proof}

So in fact the representable maps between manifolds,
in the sense of Def. \ref{def:repmap}, are exactly the submersions.
As a corollary, we see that if $p: F \ra G$ is a representable
map of stacks, the resulting map $f^{*}p$ in the diagram
\[
  \xymatrix{ M \times_{G} F \ar[r] \ar[d]_{f^{*}p} & F \ar[d]^{p} \\
             M \ar[r]_{f} & G
  }
\]
must be a submersion.

As in the topological case, a map of manifolds is an epimorphism
of stacks if and only if it has (smooth) local sections.
Hence a representable map (submersion) between manifolds is an epimorphism
of stacks if and only if it is surjective.

Much as in the topological case (and the algebraic case) we
can define many properties in terms of pullbacks, here taking care
to only use pullbacks by submersions.

\begin{Def}\label{def:propertieslocalontarget}
  Let $P$ be a property of maps of smooth manifolds that is stable under
  pullback by a submersion.
  We say that a map $f:F \ra G$ of stacks over $\bM$
  \textbf{has property $P$} if, for every representable map $g:M \ra G$ from a manifold,
  the projection $F \times_{G} M \ra M$ has property $P$.
\end{Def}
Note that this yields the same definition of property $P$ on
manifolds, since representable maps are submersions.

Hence we can define the following properties of maps of stacks: 
injective, surjective, immersion, submersion, embedding, open embedding,
closed embedding, \'etale.

We now have the smooth version of Prop. \ref{prop:presentationgivesgroupoid}.
We do not get the representability of the 
diagonal map, since that came from having arbitrary pullbacks.
\begin{proposition}\label{prop:smoothpresentationgivesgroupoid}
  Let $\GG$ be a smooth groupoid and let $\uGG$ be the 
  associated stack as in Sec. \ref{subsec:assocstack}. 
  Then $\uGG$ is locally representable; explicitly, 
  the canonical epimorphism $p:\GGo \ra \uGG$ is representable, and in particular,  
  $\GGo \times_{\uGG} \GGo \iso \GGl$. 

  Conversely, suppose that $F$ is a locally representable stack over $\bM$
  and $p:C_{0} \ra F$ is a presentation. Let $F' = C_{0} \times_{F} C_{0}$. Then
  $F'$ is represented by a smooth manifold $C_{1}$
  the pair $\mathcal{C} = (C_{1},C_{0})$ 
  has a canonical structure of a smooth groupoid,
  where the source and target maps are
  the projections $\pi_{1},\pi_{2}$; and the stack $\underline{\mathcal{C}}$
  associated to $\mathcal{C}$ is canonically equivalent to $F$.
\end{proposition}
\begin{proof}
  The proof is very similar to the proof of Prop. \ref{prop:presentationgivesgroupoid}.
  We note only the necessary changes.

  First, in the proof that $p:\GGo \ra \uGG$ is representable, we
  note that the pullbacks
  \[
    \xymatrix{ B_{\alpha} \ar[r] \ar[d] & \GGl \ar[d]^{t} \\
               U_{\alpha} \ar[r]_{a_{\alpha}} & \GGo. 
    }
  \]
  exist because $t$ is a smooth submersion. The gluing and the proof that
  the glued manifold represents the pullback are just as in the previous case;
  one simply needs to use that smoothness is local. It
  is straightforward to check that the glued manifold
  is Hausdorff and second countable.

  The proof that $\GGo \times_{\uGG} \GGo \iso \GGl$
  is identical.

  In showing the converse, we just need to note that
  the source and target maps, as pullbacks of a representable map, 
  are submersions. All of the rest of the proof is identical.
\end{proof}

\subsection{Manifolds as Sheaves}\label{subsec:manifoldassheaf}

We would like to know how to characterize manifolds
among the stacks on $\bM$, that is, how to recognize a
representable stack. Clearly any representable
stack is locally representable (use the identity as the presentation map), 
and is a sheaf. In particular the presentation map is \'etale.
These three facts are almost enough to characterize the representable
stacks, but if we want the manifold to be Hausdorff, we need an additional 
condition.

\begin{proposition}\label{prop:manifoldsassheaves}
  Let $N$ be a smooth manifold and let $F = \uN$ be the corresponding stack.
  Then
  \begin{enumerate}
  \item \label{mslr} $F$ is locally representable, with a covering $p: M \ra F$;
  \item \label{mset} the map $p$ is \'etale;
  \item \label{mssh} $F$ has trivial automorphisms (so is equivalent to a sheaf);
  \item \label{msdi} The diagonal map $\Delta: F \ra F \times F$ is a closed embedding.
  \end{enumerate}
  Conversely, if a stack $F$ on $\bM$ satisfies all of these conditions, then 
  $F$ is equivalent to a stack of the form $\uN$, for
  a smooth manifold $N$.  
\end{proposition}
\begin{proof}
  Let $N$ be a smooth manifold. We have already observed that
  $\uN$ satisfies conditions (\ref{mslr},\ref{mset},\ref{mssh}).
  Since $N$ is Hausdorff, the map $\Delta: N \ra N \times N$ is a closed embedding,
  in the usual sense. Hence it is a closed embedding in the stack sense
  (Def. \ref{def:propertieslocalontarget} and after).

  Now let $F$ be a stack satisfying the four conditions. 
  Let $M_{0} = M$ and $M_{1} = M \times_{F} M$. As in Prop.
  \ref{prop:smoothpresentationgivesgroupoid} the pair $\MM = (M_{0},M_{1})$ 
  has the structure of a groupoid. Since $F$ has trivial automorphisms,
  $(M_{0},M_{1})$ is actually an equivalence relation,
  i.e. the map $(s,t): M_{1} \ra M_{0} \times M_{0}$ is injective.
  Since the map $p$ is \'etale, so are the maps $s$ and $t$.
  Also, from Cor. \ref{cor:pullbackbydiag}, the map
  $(s,t)$ is a closed embedding.
  We know that the stack quotient of this equivalence relation
  is $F$; hence we just need to show that the stack quotient
  is the ordinary quotient, and that the ordinary quotient is
  a (Hausdorff) smooth manifold. The standard terminology
  for this is that the quotient is \textit{effective}. We will
  show this explicitly.

  Let $N$ be the usual topological quotient of $M_{0}$ by
  the equivalence relation $M_{1} \subset M_{0} \times M_{0}$.
  Since $M_{1}$ is a closed subset of $M_{0} \times M_{0}$ it
  is a standard fact that $N$ is Hausdorff. We claim
  that the quotient map $q:M_{0} \ra N$ is \'etale. For,
  we first note that $q$ is an open map, since for an 
  open set $U \subset M_{0}$,
  \[
      q^{-1}q(U) = ts^{-1}(U)
  \]
  and $t$ is an open map, so $q^{-1}q(U)$ is open and hence $q(U)$ is open.
  Also, $q$ is locally injective: given $x \in M_{0}$, there is
  a neighborhood $U$ of $x$ such that 
  $s|_{(U \times U) \intersect \Delta}$ is injective
  (since $s$ is \'etale and $M_{1}$ has the induced topology).
  But since $\Delta \subset M_{1}$ this says that 
  \[
    (U \times U) \intersect M_{1} = (U \times U) \intersect \Delta 
  \] 
  which says that $q|_{U}$ is injective. Hence $q$ is \'etale. 
  
  Now let $y \in N$ and let $x \in M_{0}$ with $q(x) = y$.
  We can find a chart $(U,\phi)$ around $x$ such that
  $q|_{U}$ is a homeomorphism onto a neighborhood $V$ of $y$.
  Hence $N$ is a topological manifold. It is easy to see that
  smooth overlap of charts in $M_{1}$ leads to smooth overlap
  in $N$. Also, since $q$ is \'etale surjective and $M_{0}$ is
  second countable, so is $N$. Hence $N$ is a smooth manifold,
  and by construction the map $q$ is \'etale in the smooth sense, i.e.
  a local diffeomorphism.
 
  Now we need to check that $N$ represents the stack quotient
  $\uMM$ of $\MM = (M_{0},M_{1})$.
  We have a natural map  $r: \uMM \ra \uN$ coming
  from $q$. We need to show that this is an isomorphism.
  First, it is an epimorphism since the map $q$ has local sections
  (being \'etale surjective).
  Second, we need to show that given any two maps $f,g: A \ra M_{0}$
  such that $qf = qg$, these maps give isomorphic maps
  to $\uMM$. But to show that we just note that since $M_{1}$ has
  the induced topology, there is a map $h:A \ra M_{1}$ such
  that $sh = f$ and $th = g$. By the definition of $\uMM$ this
  means that the maps $A \ra \uMM$ coming from $f$ and $g$ are
  isomorphic. Hence $N$ does represent the stack quotient $\uMM$
  of $\MM$, which we already know is equivalent to $F$.
\end{proof}

Of course if we do not need $N$ to be Hausdorff we do
not need to require $\Delta$ to be a \textit{closed} embedding.
This is an example of
a general phenomenon, familiar to algebraic geometers:
separation conditions on $F$ are best expressed in terms of
the diagonal map of $F$, and thus also correspond to conditions on
the map $(s,t): \GGl \ra \GGo \times \GGo$ of a representing groupoid. 

Note that in this proposition we could even require that
the manifold $M$ be an object of $\Euc$, since we can cover
any manifold by charts. 

\subsection{Orbifolds as Stacks}\label{subsec:orbifoldasstack}

Now turn to orbifolds. Following Moerdijk \cite{Moerdijk:OrbGrpIntro},
we define an orbifold as a smooth proper \'etale groupoid.
As in Section~\ref{subsec:assocstack}, we get a stack.
It is characterized by the following proposition, which
is basically just Prop. \ref{prop:manifoldsassheaves} 
with the ``trivial automorphisms'' condition removed.

\begin{proposition}\label{prop:orbifoldsasstacks}
  Let $\GG$ be an orbifold groupoid. Then the associated
  stack $F = \uGG$ satisfies the following properties:
  \begin{enumerate}
  \item $F$ is locally representable, with a covering $p: M \ra F$,
        where $M$ is a manifold;
  \item the map $p$ is \'etale;
  \item The diagonal map $\Delta: F \ra F \times F$ is proper.
  \end{enumerate}
  Conversely, let $F$ be a stack on the site $\bM$. Suppose that
  $F$ satisfies the three properties above. 
  Then $F$ is equivalent to a stack of the form $\uGG$, for
  an orbifold groupoid $\GG$.
\end{proposition}
\begin{proof}
  Let $\GG$ be an orbifold groupoid, i.e. a smooth proper     
  \'etale groupoid. Prop. \ref{prop:smoothpresentationgivesgroupoid}
  says that $F = \uGG$ is a locally representable stack.
  To show that the presentation map $p: \GGo \ra F$ is \'etale,
  we look at the construction of the pullback
  $X \times_{a,F, p} \GGo$ in 
  Prop. \ref{prop:presentationgivesgroupoid}
  and Prop. \ref{prop:smoothpresentationgivesgroupoid}.
  The pullback is constructed by gluing ordinary pullbacks
  of the map $p$ over open subsets of $X$. Since \'etale maps
  are stable under pullback and local on the target, the
  resulting map is \'etale.
 
  Similarly, the diagonal map is proper. For,
  given a \textit{representable} map $a: X \ra F \times F$,
  we can construct the pullback $X \times_{F \times F} F$
  just as in Prop. \ref{prop:presentationgivesgroupoid}.
  Properness (for locally compact Hausdorff spaces)
  is stable under pullback and local on the target
  so the resulting map is proper. Hence the diagonal is proper.  

  Now assume that we have a stack $F$ satisfying the three
  conditions above. Prop. \ref{prop:smoothpresentationgivesgroupoid}
  says that $(M,M_{1} = M \times_{F} M)$ is a smooth groupoid
  whose associated stack is equivalent to $F$. Since
  $p$ is \'etale, so are the maps $s,t:M_{1} \ra M$,
  and since $\Delta : F \ra F \times F$ is proper, so is
  $(s,t):M_{1} \ra M \times M$.
\end{proof}

In fact this is the algebraic geometer's definition of an orbifold,
as a certain kind of locally representable smooth stack.

Recall from section \ref{subsec:moritaHS} that a map of stacks associated
to groupoids corresponds to a Morita, or Hilsum-Skandalis, map
of the groupoids. Hence when we think of orbifolds as stacks
we are naturally led to the notion that a map of orbifolds
is a Morita map of the corresponding orbifold groupoids.

The distinction between an orbifold groupoid and a stack
satisfying the three conditions of Prop. \ref{prop:orbifoldsasstacks}
is roughly the same as the distinction between a particular
atlas for a manifold and the manifold itself. The groupoid
perspective is more useful for explicit computations, but the
stack perspective is more intrinsic.

\subsection{Differentiable Artin Stacks}\label{subsec:diffartin}

In the last section we saw that orbifolds correspond
to locally representable stacks with good separation
properties and an \textit{\'etale} presentation.
The analogous stacks in algebraic geometry are
known as \textit{Deligne-Mumford stacks} \cite{DeligneMumford}.
There is a more general class of stacks in algebraic
geometry known as \textit{Artin stacks}, where
one only requires a \textit{smooth} presentation.
We can make the same generalization here:

\begin{Def}\label{def:diffartin}
  A \textbf{differentiable Artin stack} is a locally
  representable stack $F$ over $\bM$ with proper
  diagonal map.
\end{Def}

Just as in the previous section, this corresponds to
a natural class of groupoids:

\begin{proposition}\label{prop:diffartin}
  A differentiable Artin stack corresponds to
  a smooth proper groupoid.
\end{proposition}
(As usual, the more precise version of this is
that there is an equivalence of bicategories, much as
in \cite{Pronk:Bicategories}.)

\begin{example}\label{ex:BGasdiffartin}
  Let $K$ be a compact Lie group. The stack
  $\BB K$ is a differentiable Artin stack,
  since it comes from the smooth proper groupoid $(\ast,K)$.
\end{example}

\begin{example}\label{ex:properactionasartin}
  Let $G$ be a Lie group acting properly on a smooth
  manifold $M$ (on the right), and let $M \rtimes G$ be the 
  translation groupoid. The associated stack $F$
  is a differentiable Artin stack.
\end{example}

\section{Stacks over Stacks}\label{sec:stackoverstack}

We now turn to some of the more subtle issues.
We want to show that the two most natural definitions of
a stack over another stack are equivalent.
However first we insert a brief discussion of a very important
kind of stack.

\subsection{Gerbes}\label{subsec:gerbes}

The stack $\BB G = \Prin_{G}$ is our primary example 
of a stack, but it is in fact a rather special stack. Note the
following elementary facts about principal bundles:
\begin{enumerate}
\item Every space $X$ has at least one principal $G$-bundle over it
      (namely, the trivial bundle).
\item Any two principal $G$-bundles are locally isomorphic.
\end{enumerate}

These two facts lead to the definition of a \textit{gerbe}.
\begin{Def}\label{def:gerbe}
  Let $(\CC,J)$ be a site and let $F:\DD \ra \CC$ be a stack.
  Then $F$ is a \textbf{gerbe} if it satisfies the following two conditions:
\begin{enumerate}
\item \label{localexist}
      For any object $C$ of $\CC$, there is a covering family
      $f_{i}:C_{i} \ra C$ such that the fiber $F_{C_{i}}$ is nonempty
      for every $i$.
\item \label{localiso}
      For any object $C$ of $\CC$ and any two objects $D_{1},D_{2}$ of
      $F_{C}$, there is a covering family
      $f_{i}:C_{i} \ra C$ such that $D_{1}|_{C_{i}}$ and $D_{2}|_{C_{i}}$
      are isomorphic.
\end{enumerate}
\end{Def}
Condition \ref{localexist} says that objects locally exist
(note this is weaker than the global existence satisfied by $\Prin_{G}$);
condition \ref{localiso} says that any two objects are locally
isomorphic.

\subsection{Equivalent definitions of a stack over a stack}
\label{subsec:equivstackstack}

Let $(\CC,J)$ be a site and consider a stack $F:\DD \ra \CC$.
We want to think of $F$ as a generalized object of
$\CC$, so we would like to know what the definition of
a stack over $F$ should be. There are two natural 
notions, and they are equivalent, as we shall now see.

The simplest definition is the following.
\begin{Def}\label{def:stackoverstack}
  Let $(\CC,J)$ be a site and consider a stack $F:\DD \ra \CC$.
  A \textbf{stack over $F$} consists of a stack $G: \EE \ra \CC$
  and a map of stacks $A: \EE \ra \DD$.
\end{Def}
So we have a commutative triangle
\[
  \xymatrix{ \EE \ar[rr]^{A} \ar[dr]_{G} &     &\DD \ar[dl]^{F} \\
                                         & \CC & 
  }
\]
Note that we have changed our usual convention of writing
a map of stacks with a lowercase letter, to put
$A$, $F$, and $G$ on an equal footing.

Alternately, one can define a Grothendieck topology $J_{\DD}$
on $\DD$ by declaring
a family of maps $\{f_{i}: D_{i} \ra D \}$ to be a covering
family if and only if $\{F(f_{i})\}$ is a covering family.
(It is easy to check that this is a topology.)
Then one could define a stack over $F$ to be a stack
$A: \EE \ra \DD$ over the induced site $(\DD,J_{\DD})$.
These notions are equivalent, by the following theorem.

\begin{theorem}\label{thm:twooutofthree}
  Let $\CC, \DD, \EE$ be three categories and let
  $F: \DD \ra \CC$, $G: \EE \ra \CC$ and $A: \EE \ra \DD$ be
  three functors with $FA = G$. Let $J$ be a Grothendieck topology
  on $\CC$ and let $J_{\DD}$ be the induced topology on $\DD$.
  Then
  \begin{enumerate}
  \item \label{compose} If $F$ is a stack over $(\CC,J)$ and $A$ is a 
        stack over $(\DD,J_{\DD})$, then $G$ is a stack over $(\CC,J)$.
  \item \label{passtoequiv} 
        If $G$ is a stack and $F$ is a prestack, then
        there is a stack $G': \EE' \ra \CC$, an equivalence
        $I: G \ra G'$ of stacks over $(\CC,J)$, and a map
        of stacks $P: G' \ra F$, with $PI=A$ such that
        the underlying functor $P:\EE' \ra \DD$ is 
        a stack over $(\DD,J_{\DD})$, as in the following diagram.
\[
\xymatrix{  
  \EE \ar[r]^{\sim}_{I} \ar[dr]_{G} \ar@/^2pc/[rr]^{A}
      & \EE' \ar[r]_{P} \ar[d]_{G'} 
      & \DD  \ar[dl]^{F} \\
    & \CC & 
  }
\]
  \item \label{stackoversheaf} 
        If $G$ is a stack over $(\CC,J)$ and $F$ is a \textit{sheaf} 
        (i.e. a discretely fibered stack) over $(\CC,J)$, then
        $A$ is a stack, i.e. we need not pass to an equivalent
        stack $G'$ as in (\ref{passtoequiv}).
  \end{enumerate}
\end{theorem}
\textit{Remark.} Since we are always willing to pass to an equivalent
stack, parts (\ref{compose}) and (\ref{passtoequiv}) above
say that the two alternate definitions of a ``stack over a stack''
are equivalent.

\begin{proof}
  \begin{enumerate}
  \item First we show that $G$ is a fibered category.
      (As always we mean ``fibered in groupoids.'')
      Let $c: C' \ra C$ in $\CC$ and let $E \in G_{C}$. 
      Then there is some $d: D \ra A(E)$ with $F(d) = c$,
      since $F$ is fibered. Since $A$ is fibered,
      there is some $e: E' \ra E$ with $A(e) = d$, hence
      $G(e) = FA(e) = c$.

      Let $e':E' \ra E$ and $e'':E'' \ra E$ in $\EE$
      and let $c: G(E') \ra G(E'')$ such that 
      $G(e'') c = G(e')$. Then there is a unique $d:A(E') \ra A(E'')$
      such that $F(d) = c$ and $A(e'') d = A(e')$, since $F$ is fibered.
      Since $A$ is fibered, there is a unique $e:E' \ra E''$
      such that $A(e) = d$ and $e'' e = e'$. Clearly
      $G(e) = c$.
      Suppose that we had another map $\tilde{e}:E' \ra E''$
      such that $G(\tilde{e}) = c$ and $e'' \tilde{e} = e'$.
      Then $FA(\tilde{e}) = c = F(d)$ and $A(e'')A(\tilde{e}) = A(e')$,
      so by the uniqueness of $d$, we must have $A(\tilde{e}) = d$.
      But then we must have $\tilde{e} = e$. Hence $G$ is fibered.

      Now we show that $G$ is a prestack. Throughout this
      proof we will avoid picking particular pullbacks,
      since that just complicates the issue. 
      Let $\{C_{\alpha} \ra C\}$ be a cover in $\CC$ and let
      $x,y \in G_{C}$. For every $\alpha$ let there be
      a diagram
\[
  \xymatrix{ x_{\alpha} \ar[r] \ar[d]        & x \\
             y_{\alpha} \ar[r] \ar[r]        & y
  }
\]
where $x_{\alpha}, y_{\alpha} \in G_{C_{\alpha}}$,
and for every pair $\alpha,\beta$ let there be a commutative diagram
as on the left:
\begin{equation}\label{prediag}
\xymatrix@=5pt{
     & & x_{\beta} \ar[dr] \ar[ddd] & 
      & & & & & & C_{\beta} \ar[dr] \ar@{=}[ddd] & \\
  x_{\alpha \beta} \ar[urr] \ar[dr] \ar[ddd] & & & x \ar@{.>}[ddd]^{\psi} 
    & & & & C_{\alpha \beta} \ar[urr] \ar[dr] \ar@{=}[ddd] & & & C \ar@{=}[ddd] \\
     & x_{\alpha} \ar[urr] \ar[ddd] & & & \ar@{~>}[rr]^{G}
       & & & & C_{\alpha} \ar[urr] \ar@{=}[ddd] & & \\
     & & y_{\beta} \ar[dr] & 
      & & & & & & C_{\beta} \ar[dr] & \\
  y_{\alpha \beta} \ar[urr] \ar[dr] & & & y 
   & & & & C_{\alpha \beta} \ar[urr] \ar[dr] & & & C \\
     & y_{\alpha} \ar[urr] & & 
     & & & & & C_{\alpha} \ar[urr] & &
  }
\end{equation}

We want to show that there is a unique $\psi: x \ra y$ in $\EE$
which fills in the dotted arrow and satisfies $G\psi = 1$.
Since $\DD$ has the induced topology, the family $\{Ax_{\alpha} \ra Ax \}$
is a cover of $Ax$, and $\{Ay_{\alpha} \ra Ay \}$ is a cover
of $Ay$. We have a commutative diagram 
\begin{equation}\label{Adiag}
\xymatrix@=5pt{
     & & A(x_{\beta}) \ar[dr] \ar[ddd] & 
      & & & & & & C_{\beta} \ar[dr] \ar@{=}[ddd] & \\
  A(x_{\alpha \beta}) \ar[urr] \ar[dr] \ar[ddd] & & & A(x) \ar@{.>}[ddd]^{\phi} 
    & & & & C_{\alpha \beta} \ar[urr] \ar[dr] \ar@{=}[ddd] & & & C \ar@{=}[ddd] \\
     & A(x_{\alpha}) \ar[urr] \ar[ddd] & & & \ar@{~>}[rr]^{F}
       & & & & C_{\alpha} \ar[urr] \ar@{=}[ddd] & & \\
     & & A(y_{\beta}) \ar[dr] & 
      & & & & & & C_{\beta} \ar[dr] & \\
  A(y_{\alpha \beta}) \ar[urr] \ar[dr] & & & A(y) 
   & & & & C_{\alpha \beta} \ar[urr] \ar[dr] & & & C \\
     & A(y_{\alpha}) \ar[urr] & & 
     & & & & & C_{\alpha} \ar[urr] & &
  }
\end{equation}
so since $F$ is a prestack, there is a unique
$\phi: A(x) \ra A(y)$ filling in (\ref{Adiag})
with $F(\phi) = 1$. Since $A$ is fibered, there is
a lift $\chi: y' \ra y$ of $\phi$ such that $A(y') = A(x)$.
We also lift the bottom square of (\ref{Adiag}) to make
a diagram 
\begin{equation}\label{prediagA}
\xymatrix@=5pt{
     & & x_{\beta} \ar[dr] \ar[ddd] & 
      & & & & & & A(x_{\beta}) \ar[dr] \ar[ddd] & & & \\
  x_{\alpha \beta} \ar[urr] \ar[dr] \ar[ddd] & & & x \ar@{.>}[ddd]_{\theta} 
    & & & & A(x_{\alpha \beta}) \ar[urr] \ar[dr] \ar[ddd] & & & A(x) \ar[ddd] & \\
     & x_{\alpha} \ar[urr] \ar[ddd] & & & \ar@{~>}[rr]^{A}
       & & & & A(x_{\alpha}) \ar[urr] \ar[ddd] & & & & \\
     & & y_{\beta}' \ar[dr] & 
      & & & & & & A(x_{\beta}) \ar[dr] & & & \\
  y_{\alpha \beta}' \ar[urr] \ar[dr] & & & y' \ar[rr]^{\chi} & & y
    & & A(x_{\alpha \beta}) \ar[urr] \ar[dr] & & & A(x) \ar[rr]^{\phi} & & A(y)\\
     & y_{\alpha}' \ar[urr] & & 
     & & & & & A(x_{\alpha}) \ar[urr] & & & &
  }
\end{equation}
Since $A$ is a prestack, there is a unique map $\theta: x \ra y'$
making the left of (\ref{prediagA}) commute and such that $A(\theta) = 1$.
Composing with $\chi$ gives our desired map $\psi: x \ra y$.
This shows existence of such a map; the uniqueness follows 
from the uniqueness of $\phi$ and the fact that $A$ is a prestack.

Now we show that $G$ is a stack. Let $\{C_{\alpha} \ra C\}$
be a cover in $\CC$. For each $\alpha$ let $x_{\alpha} \in G_{C_{\alpha}}$,
for every pair $\alpha, \beta$ let there be a diagram
\[
  \xymatrix@=5pt{     
      & x_{\alpha} & & & & & C_{\alpha} \\
    x_{\alpha \beta} \ar[ur] \ar[dr] & & \ar@{~>}[rr]^{G} & & 
        & C_{\alpha \beta} \ar[ur] \ar[dr] & \\
      & x_{\beta}  & & & & & C_{\beta} 
  }
\]
(where the arrows on the right are the canonical projections)
and for every triple $\alpha ,\beta ,\gamma$ let there be a commutative diagram
\begin{equation}\label{stackdiag}
  \xymatrix@=5pt{ 
    & x_{\alpha \beta} \ar[r] \ar[dr] & x_{\alpha} 
      & & & & & C_{\alpha \beta} \ar[r] \ar[dr] & C_{\alpha} \\
    x_{\alpha \beta \gamma} \ar[ur] \ar[r] \ar[dr] 
       & x_{\alpha \gamma} \ar[ur] \ar[dr] & x_{\beta}
      & \ar@{~>}[rr]^{G} & & & C_{\alpha \beta \gamma} \ar[ur] \ar[r] \ar[dr] 
       & C_{\alpha \gamma} \ar[ur] \ar[dr] & C_{\beta} \\
    & x_{\beta \gamma} \ar[ur] \ar[r] & x_{\gamma}
      & & & & & C_{\beta \gamma} \ar[ur] \ar[r] & C_{\gamma} 
  }
\end{equation}
where again all of the maps on the right are the canonical projections. 
We want to show that there is an $x \in G_{C}$
and arrows $\{x_{\alpha} \ra x\}$ with 
\begin{equation}\label{stackwantdiag}
  \xymatrix@=5pt{  
      & x_{\alpha} \ar[dr] & & & & & & C_{\alpha} \ar[dr] &  \\
    x_{\alpha \beta} \ar[ur] \ar[dr] & & x 
     & \ar@{~>}[rr]^{G} & & & C_{\alpha \beta} \ar[ur] \ar[dr] & & C \\
      & x_{\beta} \ar[ur]  & & & & & & C_{\beta} \ar[ur] &  
  }
\end{equation}

Applying $A$ to (\ref{stackdiag}), we get 
\begin{equation}\label{stackdiagA}
  \xymatrix@=5pt{ 
    & A(x_{\alpha \beta}) \ar[r] \ar[dr] & A(x_{\alpha}) \\
    A(x_{\alpha \beta \gamma}) \ar[ur] \ar[r] \ar[dr] 
       & A(x_{\alpha \gamma}) \ar[ur] \ar[dr] & A(x_{\beta}) \\
    & A(x_{\beta \gamma}) \ar[ur] \ar[r] & A(x_{\gamma}) \\
  }
\end{equation}
which is descent data
for $F$ on the cover $\{C_{\alpha} \ra C\}$.
Since $F$ is a stack, there is a $y \in F_{C}$
and arrows $\{A(x_{\alpha}) \ra y\}$ making the commutative diagram
\begin{equation}\label{stackwantdiagA}
  \xymatrix@=5pt{ 
      & A(x_{\alpha}) \ar[dr] & & & & & & C_{\alpha} \ar[dr] &  \\
    A(x_{\alpha \beta}) \ar[ur] \ar[dr] & & y 
       & \ar@{~>}[rr]^{F} & & & C_{\alpha \beta} \ar[ur] \ar[dr] & & C \\
      & A(x_{\beta}) \ar[ur]  & & & & & & C_{\beta} \ar[ur] &
  }
\end{equation}
Note that therefore $\{A(x_{\alpha}) \ra y\}$ is a cover, since
$\DD$ has the induced topology.

Now we claim that
$A(x_{\alpha \beta}) \iso A(x_{\alpha}) \times_{y} A(x_{\beta})$.
For, suppose we have a commutative diagram
\[
  \xymatrix@=5pt{ 
      & & A(x_{\alpha}) \ar[dr] & & & & & & & C_{\alpha} \ar[dr] & \\
    z \ar@/^1pc/[urr] \ar@/_1pc/[drr] & A(x_{\alpha \beta}) \ar[ur] \ar[dr] & & y
      & & \ar@{~>}[rr]^{F} & & F(z) \ar@/^1pc/[urr] \ar@/_1pc/[drr] 
            & C_{\alpha \beta} \ar[ur] \ar[dr] & & C \\
      & & A(x_{\beta}) \ar[ur] & & & & & & & C_{\beta} \ar[ur] & 
  }
\]
where $z$ is arbitrary.
Then there is a unique fill-in arrow $F(z) \ra C_{\alpha \beta}$,
and since $F$ is fibered, there is a unique fill-in
$z \ra A(x_{\alpha \beta})$. 
Hence $A(x_{\alpha \beta}) \iso A(x_{\alpha}) \times_{y} A(x_{\beta})$.
A similar statement applies to $A(x_{\alpha \beta \gamma})$.
Hence we can rewrite (\ref{stackdiagA}) as
\begin{equation}
\label{stackdiagAnew}
  \xymatrix{ 
      & A(x_{\alpha}) \times_{y} A(x_{\beta}) \ar[r] \ar[dr] 
        & A(x_{\alpha}) \ar[dr] & \\
    A(x_{\alpha}) \times_{y} A(x_{\beta}) \times_{y} A(x_{\gamma}) 
           \ar[ur] \ar[r] \ar[dr] 
        & A(x_{\alpha}) \times_{y} A(x_{\gamma}) \ar[ur] \ar[dr] 
        & A(x_{\beta}) \ar[r] & y \\
      & A(x_{\beta}) \times_{y} A(x_{\gamma}) \ar[ur] \ar[r] 
        & A(x_{\gamma}) \ar[ur] & \\
  }
\end{equation}
The diagram (\ref{stackdiag}) maps to this under $A$, so
it gives descent data for $A$ over the cover $\{A(x_{\alpha}) \ra y\}$.
Since $A$ is a stack, there is an $x \in A_{y} \subset G_{C}$
with 
\begin{equation}\label{stackwantdiagagain}
  \xymatrix{     & x_{\alpha} \ar[dr] &  \\
             x_{\alpha \beta} \ar[ur] \ar[dr] & & x \\
                  & x_{\beta} \ar[ur] 
  }
\end{equation}
commutative, as desired. Hence $G$ is a stack. 

\item Now we assume that $G$ is a stack and $F$ is a prestack.
      We first need to define $\EE'$ by adjoining objects
      to $\EE$; this is necessary to
      ensure that $P$ is fibered. We define
\[
  \Ob(\EE') = \{(x,a) \: | \: x \in \Ob(\EE), a: z \ra A(x) \text{~in~} \DD, 
                              \text{~with~} F(a) = 1_{G(x)} \}
\]
      and an arrow from $(x',a')$ to $(x,a)$ is simply an
      arrow of $\EE$ from $x$ to $x$.
      We note that an arrow $a$ with $F(a) = 1$ is
      necessarily an isomorphism.

      We define $I(x) = (x,1_{A(x)})$, $I(\arr{x'}{b}{x}) = b$;
      and $P(x,\arr{y}{a}{A(x)}) = y$, $P(\arr{(x',a')}{b}{(x,a)}) = c$
      where $c:y' \ra y$ is the unique arrow fitting into
      \[
        \xymatrix{ 
           y' \ar[r]^{a'} \ar@{.>}[d]_{c} & A(x') \ar[d]^{A(b)} \\
           y \ar[r]_{a}                   & A(x) 
        }
      \]
      such that $F(c) = G(b)$. 
      
      Note that $PI(x) = A(x)$ and $PI(\arr{x'}{b}{x}) = A(b)$, so $PI=A$
      as desired. The functor $I$ is clearly fully faithful (and injective 
      on objects); it also has an inverse equivalence $R:\EE' \ra \EE$
      given by $R(x,a) = x$, $R(\arr{(x',a')}{b}{(x,a)}) = b$.
      (Topologists will note that $\EE'$ is a fibrant replacement
      for $\EE$ in the model category of groupoids.)

      We need to show that $P$ is fibered. First, let $c: y' \ra y$
      and let $(x,\arr{y}{a}{A(x)}) \in \EE'$, so $P(x,a) = y$. 
      We need to find $(x',\arr{y'}{a'}{A(x')}) \in \EE'$ and 
      $b: x' \ra x$ with $P(b) = c$. But $G$ is fibered, so
      there is a $b: x' \ra x$ with $G(b) = Fc$, giving the diagram
      (solid arrows)
\begin{equation}\label{Fsquare}
\xymatrix@R=5pt{
  y' \ar@{.>}[r]^{a'} \ar[dd]_{c} & A(x') \ar[dd]^{Ab} 
    & & & & F(y') \ar@{=}[r] \ar[dd]_{F(c)} & G(x') \ar[dd]^{G(b)} \\
  & & \ar@{~>}[rr]^{F} & & & & \\
  y \ar[r]_{a} & A(x) 
    & & & & F(y) \ar@{=}[r] & G(x)
}
\end{equation}      
Since $F$ is fibered, there is a unique $a':y' \ra A(x')$ filling in the      
diagram, with $F(a') = 1$; note that $P(b) = c$ since $c$ is also unique.

Next, given a $f: x' \ra x$, $g x'' \ra x$ and a diagram of the form
\begin{equation}\label{triplediag}
\xymatrix@=5pt{ 
 y'  \ar[rr]^{a'} \ar[dr] \ar[dd]_{h} 
  & & A(x') \ar[dr]^{A(f)} \ar@{.>}'[d][dd] & 
                                        & & & &     G(x')  \ar@{=}[rr] \ar[dr] \ar[dd]_{F(h)} 
                                                          & & G(x') \ar[dr] \ar'[d]_{F(h)}[dd] & \\
  & y \ar[rr]^{a} & & A(x)       & \ar@{~>}[rr]^{F} & & &  & G(x)  \ar@{=}[rr] & & G(x)  \\
 y'' \ar[rr]_{a''} \ar[ur] 
  & & A(x'') \ar[ur]_{A(g)} &           & & & &     G(x'') \ar@{=}[rr] \ar[ur] & & G(x'') \ar[ur] & 
}
\end{equation} 
we seek to fill in the dotted arrow with an arrow of the form $A(k)$, for a unique
$k:x' \ra x''$. Since $G$ is fibered, there is a $k: x' \ra x''$ with
$gk = f$ and $G(k) = F(h)$.
Hence in the left-hand diagram of (\ref{triplediag}), the right-hand triangle
commutes, and the square in back commutes
because the dotted arrow $A(k)$ is the unique fill-in which maps to $F(h)$ under
$F$.

Hence we have found $k:(x',a') \ra (x'',a'')$ with $gk = f$ and $P(k) = h$.
It is unique, since if we have another such arrow $\tilde{k}$,
then $G (\tilde{k}) = FA(\tilde{k}) = F(h) = FA(k) = G(k)$ so $k = \tilde{k}$ by the
uniqueness of lifts through $G$.

Now that we have proved that $P$ is fibered, we can assume from now on
without loss of generality that the original functor $A$ was fibered.
This allows us to dispense with the particular construction of $\EE'$ and $P$ above,
making the argument simpler. 

We next seek to prove that $A$ is a prestack.
Let $y$ be in $\DD$ and let $\{y_{\alpha} \ra y \}$ be a cover of $y$ in the
induced topology, i.e. we have $\{F(y_{\alpha}) \ra F(y) \}$ is a cover in $\CC$.
Let $x, x'$ in $\CC_{y}$ and suppose that for each pair $\alpha ,\beta$ we have 
a diagram
\begin{equation}\label{prediag2}
\xymatrix@=5pt{
     & & x_{\beta}' \ar[dr] \ar[ddd] & 
      & & & & & & y_{\beta} \ar[dr] \ar@{=}[ddd] & \\
  x_{\alpha \beta}' \ar[urr] \ar[dr] \ar[ddd] & & & x' \ar@{.>}[ddd]^{\psi} 
    & & & & y_{\alpha \beta} \ar[urr] \ar[dr] \ar@{=}[ddd] & & & y \ar@{=}[ddd] \\
     & x_{\alpha}' \ar[urr] \ar[ddd] & & & \ar@{~>}[rr]^{A}
       & & & & y_{\alpha} \ar[urr] \ar@{=}[ddd] & & \\
     & & x_{\beta} \ar[dr] & 
      & & & & & & y_{\beta} \ar[dr] & \\
  x_{\alpha \beta} \ar[urr] \ar[dr] & & & x 
   & & & & y_{\alpha \beta} \ar[urr] \ar[dr] & & & y \\
     & x_{\alpha} \ar[urr] & & 
     & & & & & y_{\alpha} \ar[urr] & &
  }
\end{equation}
Since $G$ is a prestack and $\{F(y_{\alpha} \ra y)\}$ is a cover,
there is a unique $\psi : x' \ra x$ filling in the dotted arrow
with $G(\psi) = 1$. But then $A(\psi):y \ra y$ is the unique arrow
filling in the right-hand side of (\ref{prediag2}), so $A(\psi) = 1$
since $F$ is a prestack. Hence $\psi$ is the arrow we want,
and it is unique subject to the condition $G(\psi) = 1$, hence
certainly subject to the condition $A(\psi) = 1$.
Therefore $A$ is a prestack.

Now we show that $A$ is a stack. Let
$y$ be in $\DD$ and let $\{y_{\alpha} \ra y \}$ be a cover of $y$ as before.
Let $x_{\alpha} \in A_{y_{\alpha}}$, $x_{\alpha \beta} \in A_{y_{\alpha \beta}}$ 
and $x_{\alpha \beta \gamma} \in A_{y_{\alpha \beta \gamma}}$ and suppose
given the diagram 
\begin{equation}\label{stackdiag2}
  \xymatrix@=5pt{ 
    & x_{\alpha \beta} \ar[r] \ar[dr] & x_{\alpha} 
      & & & & & y_{\alpha \beta} \ar[r] \ar[dr] & y_{\alpha} \\
    x_{\alpha \beta \gamma} \ar[ur] \ar[r] \ar[dr] 
       & x_{\alpha \gamma} \ar[ur] \ar[dr] & x_{\beta}
      & \ar@{~>}[rr]^{A} & & & y_{\alpha \beta \gamma} \ar[ur] \ar[r] \ar[dr] 
       & y_{\alpha \gamma} \ar[ur] \ar[dr] & y_{\beta} \\
    & x_{\beta \gamma} \ar[ur] \ar[r] & x_{\gamma}
      & & & & & y_{\beta \gamma} \ar[ur] \ar[r] & y_{\gamma} 
  }
\end{equation}
Since $G$ is a stack, there is an $x$ with $G(x) = F(y)$
and for each $\alpha$, an arrow $x_{\alpha} \ra x$ with 
$G(x_{\alpha} \ra x) = F(y_{\alpha} \ra y)$.

For this $x$ we need not have $A(x) = y$. However, 
$A(x)$ fits into the same diagram in $\DD$ as $y$ does:
\begin{equation}\label{prediagnew}
\xymatrix@=5pt{
     & y_{\beta} \ar[dr] & 
      & & & & & F(y_{\beta}) \ar[dr] & \\
  y_{\alpha \beta} \ar[ur] \ar[dr] & & y 
    & \ar@{~>}[rr]^{F} & & & F(y_{\alpha \beta}) \ar[ur] \ar[dr] & & F(y) \\
     & y_{\alpha} \ar[ur] & & 
       & & & & F(y_{\alpha}) \ar[ur] & 
  }
\end{equation}
Hence there is a unique isomorphism $\chi: A(x) \ra y$
which is compatible with the identity maps on the $y_{\alpha}$
and which satisfies $F(\chi) = 1$. 
Since we are assuming $A$ is fibered, the map $\chi$ lifts to
an isomorphism $\theta: x \ra x'$ for some $x'$. Using
$x'$ instead of $x$, we have our desired amalgamation of the $x_{\alpha}$'s
satisfying $A(x') = y$.

\item If $F$ is a sheaf, so $\DD$ is discretely fibered over
$\CC$, then the category $\EE'$ constructed in the previous part
is isomorphic to $\EE$, since any arrow $a$ in $\DD$ with $F(a) = 1$ is
necessarily an identity map. Hence $A$ is already fibered,
with no replacement necessary.

\end{enumerate}
\end{proof}

In particular, we have the following corollary (which can be proved
quite simply on its own) about stacks over objects of $\CC$.

\begin{cor}\label{cor:stackoverobject}
  Let $(\CC,J)$ be a site and let $C$ be an object of $\CC$.
  Then the following data are equivalent:
  \begin{enumerate}
  \item a stack $F$ over $(\CC,J)$ and a map of stacks $a:F \ra \uC$.
  \item a stack $G$ over the site $(\CC/C,J/C)$ where
        $J/C$ is the induced topology on $C$.
  \end{enumerate}
\end{cor}
Hence these two alternate notions of ``a stack over the object $C$''
are equivalent. However note that in the topological case,
there is another notion of a stack over a space $X$, namely,
a stack over the \textit{small} site $\Op(X)$. This
is different from the above two (large) notions, although it
is closely related. See Section~\ref{sec:largesmall} below.

Now we can give an alternate definition of a gerbe.
We make it in a seemingly more general context than the
previous definition.

\begin{Def}\label{def:gerbeoverstack}
  Let $(\CC,J)$ be a site and let $F: \DD \ra \CC$ be a stack.
  A \textbf{gerbe over} $F$ is a stack $G: \EE \ra \CC$ over $\CC$
  and a map of stacks $a: G \ra F$ such that:
  \begin{enumerate}
  \item \label{objepi} $a$ is an epimorphism;
  \item \label{arrepi} 
        the diagonal map $\Delta_{a}: G \ra G \times_{F} G$ is an epimorphism.
  \end{enumerate}
\end{Def}

This definition is in fact equivalent to Definition~\ref{def:gerbe}.
\begin{proposition}\label{prop:gerbedefsequiv}
  Let $(\CC,J)$ be a site and let $F: \DD \ra \CC$ be a stack.
  The data of a gerbe over $F$ (Def. \ref{def:gerbeoverstack}) 
  is equivalent to the data of a
  gerbe over the induced site $(\DD,J_{\DD})$ (Def. \ref{def:gerbe}).
  In particular, a gerbe over $(\CC,J)$ in the sense of Def. \ref{def:gerbe}
  is equivalent to a gerbe over the final stack $F = \mathrm{Id}_{\CC}$
  on $\CC$.
\end{proposition}
\begin{proof}
Thm. \ref{thm:twooutofthree} ensures that the data of
a stack over $F$ and a stack over $(\DD,J_{\DD})$ are equivalent.
We just need to check that the gerbe conditions correspond.
First, suppose that $G$ satisfies (\ref{objepi}) of Def. \ref{def:gerbeoverstack}.
Then given $y \in F_{C} \subset \DD$, there is some cover $\{C_{\alpha} \ra C\}$
in $\CC$ and objects $x_{\alpha} \in \EE$ with $a(x_{\alpha}) \iso y|_{C_{\alpha}}$.
But $\{y|_{C_{\alpha}} \ra y \}$ is a cover of $y$ for the induced topology,
so $a$ satisfies (\ref{localexist}) of Def. \ref{def:gerbe}.
Conversely, suppose that $G$ satisfies (\ref{localexist}) of Def. \ref{def:gerbe}.
Let $y \in F_{C}$. Then there is a cover $\{y_{\alpha} \ra y\}$ in the
induced topology on $\DD$ and objects $x_{\alpha} \in a_{y_{\alpha}}$
with $a(x_{\alpha}) = y_{\alpha}$. Since $F(\{y_{\alpha} \ra y \})$ is
a cover in $C$, this says that $a$ is an epimorphism.

Next, suppose that $G$ satisfies (\ref{arrepi}) of Def. \ref{def:gerbeoverstack}.
Let $y \in F_{C}$ and let $x,x' \in a_{y}$. 
Recall that an object of $G \times_{F} G$ over $C$ is a triple
$(x,x',\phi)$ with $x, x' \in G_{C}$ and $\phi: a(x) \ra a(x')$ in $F_{C}$.
Hence $(x,x',1_{y})$ is an object of $G \times_{F} G$ over $C$.
Since $\Delta : G \ra G \times_{F} G$
is an epimorphism, there is a cover $\{C_{\alpha} \ra C\}$ and 
objects $z_{\alpha} \in G_{C_{\alpha}}$ with 
$(z_{\alpha},z_{\alpha},1) \iso (x|_{C_{\alpha}},x'|_{C_{\alpha}},1)$,
which implies $x|_{C_{\alpha}} \iso x'|_{C_{\alpha}}$. Hence
$a$ satisfies (\ref{localiso}) of Def. \ref{def:gerbe}.
Again the converse is straightforward, and we leave it to the reader.
\end{proof}

In particular, we know what a gerbe over a sheaf is. In fact,
any stack can be expressed as a gerbe over a sheaf.
For, given a stack $F$ over $\CC$, define the \textit{sheaf quotient}
$\EG(F)$ to be the sheafification of the presheaf
\[
  \PEG(C) = \{\text{isomorphism classes of objects of~} F_{C} \}).
\]
(The notation $\EG(F)$ is from the French term 
\textit{espace grossier}, which translates as ``coarse quotient.''
However we do not wish to confuse this quotient with other
notions of coarse quotient.) 
Then the natural map $F \ra \EG(F)$
defines a gerbe. Also, it is universal for maps of
stacks to sheaves (and hence in particular to representable stacks).
\begin{proposition}\label{prop:sheafquotient} 
  Let $F$ be a stack over $\CC$, let $\EG(F)$ be its
  sheaf quotient and let $p: F \ra \EG(F)$ be the natural map.
  Then the map $p$ defines $F$ as a gerbe over $\EG(F)$.
  Given any sheaf $G$ over $\CC$ and a map of stacks $f:F \ra G$,
  $f$ factors uniquely through $p$.
\end{proposition}
\begin{proof}
  Clearly $p$ makes $F$ into a gerbe over $\EG(F)$,
  since locally, any $x \in \EG(F)(C)$ is represented
  by an object of $F$, and locally, any two such
  representatives are isomorphic, by construction.

  Given a sheaf $G$ over $\CC$ and a map of stacks $f:F \ra G$,   
  define $g: \EG(F) \ra G$ as follows. There is
  an obvious map of presheaves $\tilde{g}: \PEG(F) \ra G$ defined by
  $\tilde{g}([x]) = f(x)$. Let $g$ be the
  map of sheaves induced by $\tilde{g}$.
  Then we have $f = gp$. The uniqueness of $\tilde{g}$ is
  clear, hence so is the uniqueness of $g$.
\end{proof}

\section{Large and Small Sheaves and Stacks}
\label{sec:largesmall}

Fix a space $X$ in $\bT$. 
We would like to make clear the relationship between
stacks on the large topological site $\bTX$ and on
the small site $\Op(X)$. The former should be considered
as generalized spaces over $X$, the latter as 
algebraic objects over $X$. 

\subsection{Sheaves}\label{subsec:largesmallsheaves}

First, let us address the case of sheaves, which is better-known
and simpler. We seek to compare the categories $\Sh(\bTX)$
and $\Sh(\Op(X))$ of sheaves (of sets). There is an obvious restriction
functor $R: \Sh(\bT) \ra \Sh(\Op(X))$ coming from the inclusion
$\Op(X) \ra \bTX$.

Going the other way, given a sheaf $F$ over $\Op(X)$
and a map $f:Y \ra X$, we get a pullback sheaf $f^{*}F$ over
$\Op(Y)$ in the usual way. Define $PF(Y \stackrel{f}{\ra}) = f^{*}F(Y)$,
the global sections of the pulled-back sheaf. Given another map
$f': Y' \ra X$ and a map $g:Y' \ra Y$ with $fg=f'$,
we have a natural map of sheaves over $Y$
\[
  f^{*} F \ra g_{*} g^{*} f^{*} F \iso g_{*} (f')^{*} F
\]
and, taking global sections on $Y$, we get a map
\[
  PF(Y \stackrel{f}{\ra} X) \ra PF(Y' \stackrel{f'}{\ra} X).
\]
It is easy to see that this defines a presheaf over $\bTX$.
In fact it is a sheaf. Also, it is easy to define the action of
of $P$ on maps of sheaves. This defines a ``prolongation''
functor $P: \Sh(\Op(X)) \ra \Sh(\bTX)$.

\begin{lemma}[\cite{SGA4.1}]\label{lem:prolongation}
  Let $X$ be a space in $\bT$ and let $R: \Sh(\bTX) \ra \Sh(\Op(X))$
  and $P: \Sh(\Op(X)) \ra \Sh(\bTX)$ be the restriction and
  prolongation functors defined above. Then $P$ is left adjoint to
  $R$, and both functors are both left and right exact.
  The functor $P$ is fully faithful and the adjunction map 
  $1 \ra RP$ is an isomorphism. 
\end{lemma}

There is another, more topological, description of the
prolongation functor. Recall that to any sheaf 
$F$ on $\Op(X)$ we can associate the \textit{\'etale space}
$E(F)$, see e.g. \cite{MacMoer:Sheaves}. This is a space $\pi: E(F) \ra X$ over $X$,
so it represents a sheaf over $\bTX$. Explicitly, 
\[
  \underline{E(F)}(\arr{Y}{f}{X}) = \Maps_{X}(Y,E(F)).
\]
This can also be expressed as
\[
  \Maps_{X}(Y,E(F)) = \Gamma(Y,f^{*}E(F)) = \Gamma(Y,f^{*}F)
\]
(since the pullback of the sheaf agrees with the pullback
of the \'etale space). Hence this agrees with the prolongation
functor defined above:
\[
  P(F) = \underline{E(F)}.
\]

From the latter description we see that the image of $P$
is the full subcategory of large sheaves over $X$ which
are representable, and whose representing space is \'etale
over $X$. So from the point of view of generalized spaces,
small sheaves are quite special. The extreme example is when
$X = \ast$, a one-point space. Then $\Sh(\Op(X)) = \Sets$,
whereas $\Sh(\bT/\ast) = \Sh(\bT)$, the category of large
absolute sheaves, which includes the category of all spaces
(and much more). So in this sense the category of small sheaves
is a rather meager subcategory of the category of large sheaves.

However, from the perspective of sheaves as algebraic objects
over $X$, the two categories are quite similar. This is because
the categories $\Sh(\Op(X))$ and $\Sh(\bTX)$, while not
equivalent, are \textit{homotopy equivalent} in the sense
of topos theory. We briefly outline the situation, as presented in 
\cite{SGA4.1}.

We noted above that $R$ is both left and right exact. In fact
it is not hard to see that it preserves all limits and all
colimits. By general principles \cite{SGA4.1}, it must have a
right adjoint $Q$ as well as the left adjoint $P$.
(Unfortunately it is hard to write down $Q$ explicitly, see \cite{MacMoer:Sheaves}.)
We then have
\begin{proposition}\label{prop:cohomologyequivalent}
  Let $F$ be a sheaf of abelian groups over $\bTX$. Then there is 
        a canonical isomorphism of cohomology groups 
        \begin{equation}\label{Rcoho}
          H^{i}(\bTX,F) \rai H^{i}(X,RF).
        \end{equation}

  Let $G$ be a sheaf of abelian groups over $\Op(X)$. 
        Then there is a canonical isomorphism of cohomology groups 
        \begin{equation}\label{Pcoho}
          H^{i}(X,G) \rai H^{i}(\bTX,PF).
        \end{equation}
\end{proposition}
\begin{proof}

  The restriction functor $R: \Sh(\bTX) \ra \Sh(\Op(X))$
        is exact and takes injectives to injectives (since BLAH).
        Therefore its right derived functors vanish, and $R$ induces an
        isomorphism
        \begin{equation}\label{cohoequal}
          H^{i}(\bTX,F) \rai H^{i}(X,RF).
        \end{equation}

  To prove (\ref{Pcoho}), 
  apply (\ref{cohoequal}) to the functor $PG$ and use the
  canonical isomorphism $RPG \iso G$.
\end{proof}

So if our intention is to analyze $X$ by using cohomology of sheaves,
we can use either small sheaves or large sheaves. But
if we interested primarily in the sheaves themselves,
there is a big difference betweem small and large sheaves.

\subsection{Stacks}\label{subsec:largesmallstacks}

Now we turn to stacks. We will see that as in the sheaf case,
small stacks form a fairly special subcategory of large stacks.

Again we fix a space $X$ and consider the two 2-categories
$\St(\bTX)$, $\St(\Op(X))$ of large and small stacks respectively.
First, we still have an obvious restriction functor
$R: \St(\bTX) \ra \St(\Op(X))$.

We would like to construct a natural functor $P$ going the other way
such that $RP \iso 1$. There are two ways to do this (as
in the sheaf case). For concreteness we will emulate the
\'etale space construction given above.

First, we note without proof a result of Giraud about strictifying
stacks. 
\begin{proposition}[Giraud \cite{Giraud:Book}]
\label{prop:strictify}
  Let $F$ be a stack over a site $\CC$. Then there is a   
  strict stack $G$ over $\CC$ and a natural equivalence $i:F \rai G$
  of stacks.
\end{proposition}
In other words every stack is naturally equivalent to
a strict stack. So if necessary, when doing a construction
on stacks, we can assume that the stacks in question are strict.
However it should be noted that the strictification is often
more complicated to work with in detail than the original stack.

Hence we can suppose that we start with a strict small stack $F$ on 
$\Op(X)$. In particular this is an honest sheaf of groupoids.
It is easy to see (\cite{SGA4.1}, \cite{MacMoer:Sheaves})
that this is the same thing as a groupoid object in the category
$\Sh(\Op(X))$ of (small) sheaves of sets on $X$. 
But this is the same as a groupoid object, call it $\GG$, in the
category of \'etale spaces over $X$. Explicitly, we have
\[
  \GGo = E(\Ob F), \qquad  \GGl = E(\Ar F).
\]
Note that we have a commutative diagram 
\[
 \xymatrix{ 
 \GGl \ar[rr]^{s} \ar[dr]_{\mathrm{pr}_{1}} & &\GG0 \ar[dl]^{\mathrm{pr}_{0}}\\
                                           & X & 
 }
\]
where $\mathrm{pr}_{1},\mathrm{pr}_{2}$ are \'etale; hence
$s$ is also \'etale, so $\GG$ is an \'etale groupoid.

\begin{example}\label{ex:smallprincipalgerbe}
Let $K$ be a topological group and let 
$F$ be the small gerbe of $K$-principal bundles
over $X$. Note that this is already a strict stack (since
restriction to open subsets of $X$ commutes on the nose). 
So it is a sheaf of groupoids, with 
\[
  F(U) = \{\text{groupoid of~}K-\text{bundles over~} U \}.
\]
Expressed as a groupoid object in $\Sh(\Op(X))$,
this is 
\begin{align*}
  F_{0}(U) &= \{K-\text{bundles over~} U \} \\
  F_{1}(U) &= \{\text{maps of~} K-\text{bundles over~} U \}.
\end{align*}
\end{example}

We have so far defined a functor
\[
  P_{1}: (\text{Strict small stacks over~} X) \ra 
           (\text{\'Etale groupoids, \'etale over~} X).
\]
In fact this functor is an embedding (it just forgets the fact
that $F$ is a stack).

Now we observe that an \'etale groupoid $\GG$ which is 
\'etale over a space $X$ is quite special.
First we define the local action of a stabilizer arrow
in an \'etale groupoid.

The \textit{stabilizer space} $S_{\GG}$ of a groupoid $\GG$ 
(denoted simply by $S$ when it does not cause confusion) is
defined to be the set of arrows with the same source and target:
\[
  S_{\GG} = \{ g \in \GG_{1} \: | \: s(g) = t(g) \}.
\]
(In categorical language, these are the automorphisms.)
There is a right action of $\GG$ by conjugation: given $g \in \GG_{1}$,
$h \in S_{\GG}$, then we define $h \cdot g = g^{-1}hg$.
Hence $S_{\GG}$ is a $\GG$-space.
Given $x \in \GGo$, let $S_{x} = \{g:x \ra x \}$ be the space of
stabilizer arrows at $x$; if $\GG$ is \'etale then $S_{x}$ is discrete.

Let $\GG$ be an \'etale groupoid. An arrow $g: x \ra y$ induces
a germ of a homeomorphism $\phi_{g}$ from $x$ to $y$. 
For, choose an open neighborhood $U$ of $g$ such that
$s|_{U}: U \ra s(U) \subset \GG_{0}$ is a homeomorphism.
Then we have a continuous map $t \circ (s|_{U})^{-1}: s(U) \ra t(U)$,
such that $x \mapsto y$.
Different choices of $U$ give the same germ. This defines a 
strict homomorphism $\phi:\GG \ra \mathrm{Homeo}(\GGo)$,
the latter being the groupoid of germs of homeomorphisms of $\GGo$.

In particular, given a stabilizer arrow $(g:x \ra x) \in S$, consider the
germ $\phi_{g}$, which is the germ of a diffeomorphism fixing $x$.
The space of \textit{ineffective stabilizers} $S^{0} = S^{0}_{\GG}$ 
is defined as the subspace of $S$ consisting of the arrows 
which induce the trivial germ (the germ of the identity map)
on $\GG_{0}$. 
Let $S^{0}_{x} = S^{0}_{\GG} \intersect S_{x}$ be the set of 
ineffective stabilizers at $x$.

\begin{proposition}\label{prop:PI}
  Let $\GG$ be a groupoid in the category of \'etale spaces over $X$.
  Then any stabilizer arrow in $\GG$ must act ineffectively.
\end{proposition}
\begin{proof}
  Let $\pi: \GGo \ra X$ and $\tilde{\pi}: \GGl \ra X$ be the \'etale projections to $X$. 
  Let $x \in \GGo$ and let $g: x \ra x$. Let $U \subset \GGo$ be a neighborhood of $x$ 
  such that $\pi|_{U}:U \ra \pi(U)$ is a homeomorphism. 
  In particular $\pi$ is injective. Since
  the structure maps of $\GG$ must respect the projections, 
  we have $\pi s = \tilde{\pi} = \pi t$. Hence for any arrow $h$
  with source and target in $U$, we have $s(h) = t(h)$
  since $\pi|_{U}$ is injective. Therefore the
  germ $\phi_{g}$ is the identity.
\end{proof}
We call an \'etale groupoid where every arrow acts ineffectively
a \textit{purely ineffective} groupoid. For such a
groupoid, we have the following.
\begin{lemma}[\cite{HenriquesMetzler:OrbiPres}]\label{lem:PIquotientOK}
  Let $\GG$ be a purely ineffective \'etale groupoid. Then
  the corresponding effective groupoid $\Geff$ is
  equivalent to the topological quotient $\Gtop$
  and the projection map $p:\GGo \ra \Gtop$ is \'etale.
\end{lemma}

In our current situation, we have the following diagram:
\[
  \xymatrix@1{ 
    \GGl \ar@<1ex>[r]^{s} \ar@<-1ex>[r]_{t} 
       & \GGo \ar[r]^{p} & \Gtop \ar[r]^{q} & X.
  }
\]
The composition $qp$ of the last two maps is the \'etale
map $\mathrm{pr}_{0}$, so $q$ is also \'etale.

Hence we have two pieces of data: a purely ineffective \'etale groupoid $\GG$
over its topological quotient $\Gtop$; and an \'etale map
$q: \Gtop \ra X$. The latter is just a small sheaf over $X$. (Note
however that $q$ may not be surjective, so some stalks
can be empty.) So we will concentrate on the groupoid.
As a groupoid over $\Gtop$, it represents a large stack $\uGG$
over $\Gtop$ as in Section~\ref{subsec:assocstack}.
Since $\GG$ is purely ineffective and $\Gtop$ is its quotient,
we actually get more.
\begin{proposition}\label{prop:PIgivesgerbe}
  Let $\GG$ be a purely ineffective \'etale
  groupoid and let $\Gtop$ be its topological quotient. 
  The large stack $\uGG$ over $\Gtop$ is a locally representable gerbe with
  discrete stabilizers. 
\end{proposition}
\textit{Note.} 
In general, we say that a locally representable large 
stack $F$ \textit{has discrete stabilizers} if there is a presentation of
$F$ by a groupoid with discrete stabilizers. It easily follows that
every presentation of $F$ must have discrete stabilizers.
\begin{proof}
  Since $\uGG$ comes from the groupoid $\GG$ it is clearly locally
  representable, and since $\GG$ is \'etale, $\GG$, and hence $\uGG$,
  have finite stabilizers.

  By Prop.~\ref{lem:PIquotientOK},
  $\GG_{0} \ra \Geff = \Gtop$ always has local sections (being \'etale),
  verifying condition (\ref{localexist}) of Def. \ref{def:gerbe}.
  Further, the map $\GG_{1} \ra \GG_{0}$ has local sections,
  and hence any two objects of $\uGG(U)$ are locally isomorphic,
  verifying condition (\ref{localiso}).
  Therefore $\uGG$ is indeed a gerbe. 
\end{proof}
Note that if there are effective
stabilizers, one will not get local sections of $\GG_{0} \ra \Gtop$,
so one does not obtain a gerbe over $\Gtop$.

In fact we have now identified the large stacks which come from
small stacks in this way. First, we note the following.
\begin{theorem}\label{thm:whichlargegerbe}
  Given a space $Y$, the following three categories 
  are equivalent:
  \begin{enumerate}
  \item The category of small gerbes over $Y$;
  \item The category of purely ineffective \'etale topological groupoids 
        with $\Gtop \iso Y$;
  \item The category of locally representable large gerbes over $Y$ with
        discrete stabilizers.
  \end{enumerate}
\end{theorem}
\begin{proof}
  We will concentrate on the objects and leave the rest of
  the verification of the equivalences to the reader.
  In fact these are more precisely equivalences of bicategories.

  Let $F$ be a small gerbe over $Y$, asssumed without
  loss of generality to be strict. As we saw above, this
  determines a groupoid in the category of \'etale spaces
  over $Y$; call it $\GG$. The projections $\GGl \ra Y$ and
  $\GGo \ra Y$ are \'etale, so the maps $s,t: \GGl \ra \GGo$
  are \'etale. Since $F$ is a gerbe, the projection
  $\GGo \ra Y$ is surjective, and the natural map 
  \[
    \GGl \ra \GGo \times_{Y} \GGo
  \]
  is surjective. This implies that the natural map $\Gtop \ra Y$
  is an isomorphism. As observed above,
  since the quotient map $\GGo \ra \Gtop$ has local sections,
  $\GG$ must be purely ineffective. 

  Conversely, given a purely ineffective \'etale groupoid $\GG$, we have
  seen that this a groupoid object in the category
  of \'etale spaces over $\Gtop \iso \Geff$, in other words,
  a (strict) small stack. It is a small gerbe by the same argument
  as in Prop. \ref{prop:PIgivesgerbe}.

  Such a groupoid $\GG$ also 
  determines a locally representable large gerbe
  with discrete stabilizers, by Prop. \ref{prop:PIgivesgerbe}.
  Last, let $F$ be a large gerbe over $Y$. By restricting to the
  subcategory $\Op(Y)$, this determines a small gerbe over $Y$.
  This is clearly compatible with representing $F$ by
  a groupoid $\GG$ over $Y$.
\end{proof}

Now we include the map $q: \Gtop \ra X$ and get the following
proposition. Let $P$ be the functor which associates
to a small (strict) stack $F$ over $X$ the large stack
represented by the \'etale groupoid $\GG = P_{1}(F)$.
Putting Thm. \ref{thm:whichlargegerbe} together with
the sheaf $Y \ra X$ gives the following characterization
of the large stacks over $X$ which come from small stacks
over $X$:
\begin{theorem}\label{thm:whichlarge}
  Given a space $X$, a large stack $F$ over $X$ arises
  as the prolongation of a small stack over $X$
  if and only if it satisfies the 
  following four conditions:
\begin{enumerate}
\item $F$ is locally representable;
\item $F$ has discrete stabilizers;
\item the sheaf quotient $\EG(F)$ is representable;
\item the induced map $EG(F) \ra X$ (Prop.~\ref{prop:sheafquotient}) 
      is \'etale.
\end{enumerate}
  The large stack $F$ corresponds to a small \textit{gerbe} over
  $X$ if and only if, in addition to these conditions, 
  the map $EG(F) \ra F$ is an isomorphism. 
\end{theorem}

We now present two examples to show the
necessity of the conditions in the proposition.
\begin{example}\label{ex:largegerbenotsmallgerbe}
Consider the large absolute gerbe $\BB K$, where $K$ is a nondiscrete
topological group. This is a large gerbe over $\ast$.
However a small gerbe over $\ast$ is just an abstract group, with
no topology. So there is no way to recover $\BB K$ from
its restriction to $\Op(\ast)$. 
\end{example}

\begin{example}\label{ex:largesheafnotetale}
Let $\uZ$ be the large \textit{sheaf} over a space $Y$
represented by the space $f:Z \ra Y$ over $Y$. Since it is a
sheaf, it trivially has discrete stabilizer. However we know
it can only come from a small sheaf over $Y$ if the
map $f$ is \'etale. 
\end{example}

We should mention that one can define the equivalence
of Thm.~\ref{thm:whichlarge} by a process analogous
to the first construction in Sec.~\ref{subsec:largesmallsheaves}.
Given a continuous map $f: Y \ra X$, one can define
the pullback $f^{*}F$ of a small stack $F$ on $X$.
If we do this for all such maps $f$ and take global
sections on each $Y$, we get a large stack over $X$,
much as in Sec.~\ref{subsec:largesmallsheaves}.
We get a functor that is equivalent
to the functor $P$ defined above using \'etale spaces.



\end{document}